\documentclass[12pt,draftcls,onecolumn]{IEEEtran}
\usepackage{subfigure,amsmath,amsfonts,amssymb,color,graphicx, algorithmic,verbatim,mathrsfs,multirow,cite}

\newtheorem{lemma}{Lemma}

\newtheorem{proposition}{Proposition}

\newtheorem{assumption}{Assumption}

\newcommand{\qed}{\nobreak \ifvmode \relax \else
      \ifdim\lastskip<1.5em \hskip-\lastskip
      \hskip1.5em plus0em minus0.5em \fi \nobreak
      \vrule height0.75em width0.5em depth0.25em\fi}

\def\dist{\mathrm{dist}}
\def\ave{\mathrm{ave}}
\def\proj{\mathsf{\Pi}}

\def\a{\alpha}
\def\b{\beta}
\def\e{\epsilon}
\def\g{\gamma}
\def\t{\theta}

\def\Xc{\mathcal{X}}
\def\Yc{\mathcal{Y}}
\def\Fc{\mathcal{F}}
\def\Nc{\mathcal{N}}
\def\Es{\mathsf{E}}
\def\Ps{\mathsf{\Pi}}
\def\1b{\mathbf{1}}
\def\EXP#1{{\mathsf{E}\left[#1 \right]}}
\def\la{\langle}
\def\ra{\rangle}

\begin{document}
\title{\bf Asynchronous Gossip-Based Random Projection Algorithms Over Networks}
\author{Soomin Lee and Angelia Nedi\'{c}\thanks{
S.~Lee is with the Electrical and Computer Engineering Dept.\ and A.~Nedi\'c is with
the Industrial and Enterprise Systems Engineering Dept.\ at
the University of Illinois at Urbana-Champaign, Urbana, IL 61801;
e-mails: \text{\{lee203, angelia\}}@illinois.edu.
The support for this work has been provided by the NSF under grants
CMMI 07-42538 and CCF 11-11342.
}}

\maketitle
\begin{abstract}
We consider a fully distributed constrained convex optimization problem over a multi-agent (no central coordinator) network.
We propose an asynchronous gossip-based random projection (GRP) algorithm that solves the distributed problem
using only local communications and computations.
We analyze the convergence properties of the algorithm for an uncoordinated diminishing stepsize
and a constant stepsize.
For a diminishing stepsize, we prove that
the iterates of all agents converge to the same optimal point with probability 1.
For a constant stepsize, we establish an error bound on the expected distance from
the iterates of the algorithm to the optimal point.
We also provide simulation results on a distributed robust model predictive control problem.
\end{abstract}

\section{Introduction}
A number of important problems that arise in various application domains,
including distributed control \cite{Camponogara}, large-scale machine learning \cite{Langford-book,Bianchi-book}, wired and wireless networks \cite{johansson,ram_info,rabbat,con03}
can be formulated as a distributed convex constrained minimization problem over a multi-agent network.
The problem is usually defined as a sum of convex objective functions over an intersection of convex constraint sets.
The goal of the agents is to solve the problem in a distributed way, with each agent handling a component of the objective and constraint.
This is useful either when the problem data are naturally distributed
or when the data are too large to be conveniently processed by a single agent.
Common to these distributed optimization problems are the following operational restrictions:
1) a component objective function and a constraint set is only known to a specific network agent (the problem is fully distributed);
2) there is no central coordinator that synchronizes actions on the network or works with global information;
3) the agents usually have a limited memory, computational power and energy;
and 4) communication overhead is significant due to the expensive start-up cost and network latencies.
These restrictions motivate the design of distributed, asynchronous, computationally simple and local communication based algorithms.

The focus of this paper is the development and analysis of an efficient distributed algorithm
whereby only a pair of agents exchanges local information and updates in an asynchronous manner.
We propose a gradient descent with {\it random projections} which uses gossip scheme as a communication protocol.
Random projection-based algorithms
have been proposed in~\cite{SL12b} (see also its extended version~\cite{Soomin-report})
for distributed problems with a synchronous update rule, and in~\cite{AN2011} for centralized problems.
Synchronous algorithms are often inefficient as they create
bottlenecks and waste CPU cycles, while centralized approaches are inapplicable in situations where a central coordinator does not exist.
Asynchronous algorithms based on a gossip scheme
have been proposed and analyzed for a scalar objective function and a diminishing stepsize~\cite{RamNV09},
and a vector objective function and a constant stepsize~\cite{RamNV10}.
An asynchronous broadcast-based algorithm has also been proposed in~\cite{Nedic11}.
The gradient-projection algorithms proposed in the papers~\cite{Nedic2009, RamNV09,RamNV10,Nedic11,Bianchi2012,tsianos2012consensus}
assume that the agents share a common constraint set
and the projection is performed on the whole constraint set at each iteration.
To accommodate the situations where the agents have local constraint sets, the distributed gradient methods
with distributed projections on local constraint sets have been considered in~\cite{Nedic2010,KS2011}
(see also~\cite{Kunal-thesis}). However, even
the projection on the entire (local) constraint set often overburdens agents, such as wireless sensors,
as it requires intensive computations. Furthermore,
in some situations, the constraint set can be revealed only component-wise in time, and the whole set is not available in advance, which makes the existing distributed methods inadequate. Our proposed algorithm is intended to accommodate such situations.

In our algorithm, we efficiently handle the projection at each iteration by
performing a projection step on the local constraint set that is randomly selected (by nature or by an agent itself).
For asynchrony, each agent uses either a diminishing stepsize that is uncoordinated with those of the other agents or a constant stepsize.
Our main goals are to establish the convergence of the method with a diminishing stepsize, to estimate the error bound for a constant stepsize, and to provide simulation results for the algorithm.

To the best of our knowledge, there is no previous work on asynchronous distributed optimization algorithms that
utilize random projections.
Finding probabilistic feasible solutions through random sampling of constraints for optimization problems with uncertain constraints have been proposed in \cite{tempo2009,Calafiore:2010}.
Also, the related work is the (centralized) random projection method proposed
by Polyak \cite{Polyak2001409} for a class of convex feasibility problems and the random projection algorithm \cite{nedicFP2011} for convex set intersection problems.
On a broader scale, the work in this paper is related to the literature on the consensus problem
(see for example\cite{con01,con03,con05}).

The rest of this paper is organized as follows.
In Section \ref{sec:algo}, we describe the problem of interest, propose our gossip-based random projection algorithm, and state assumptions on the problem and the network.
Section \ref{sec:main} states the main results of the paper, while
in Sections~\ref{sec:conv} and~\ref{sec:err}, we provide the proofs
of the results.
We present the simulation results on a distributed model predictive control problem in
Section~\ref{sec:sim}
and conclude with a summary in Section \ref{sec:con}.
Appendix contains the proofs of the lemmas given in Section \ref{sec:conv} and Section \ref{sec:err}.

\noindent\textbf{Notation}.
A vector is viewed as a column. We write $x'$ to denote the transpose of a vector $x$.
The scalar product of two vectors $x$ and $y$ is $\langle x, y \rangle$.
We use $\1b$ to denote a vector whose entries are 1 and $\|x\|$ to denote the standard Euclidean norm.
We write $\dist(x,\mathcal{X})$
for the distance of a vector $x$ from a closed convex set $\mathcal{X}$, \textit{i.e.},
$\dist(x,\mathcal{X}) = \min_{v\in\mathcal{X}}\|v-x\|$. We use $\mathsf{\Pi}_{\mathcal{X}}[x]$
for the projection of a vector $x$ on the set $\mathcal{X}$, \textit{i.e.}, $\mathsf{\Pi}_{\mathcal{X}}[x] = \arg\min_{v\in\mathcal{X}}\|v-x\|^2$.
We use $\mathsf{E}[Z]$ to denote the expectation of a random variable $Z$.
We often abbreviate \textit{with probability 1} as \textit{w.p.1}.

\section{Problem Set-up, Algorithm and Assumptions \label{sec:algo}}
We consider an optimization problem where the objective function and constraint sets are distributed among $m$ agents over a network.
Let an undirected graph $G=(V,E)$ represent the topology of the network, with the vertex set $V = \{1,\ldots,m\}$ and the edge set $E \subseteq V \times V$. Let $\Nc(i)$ be the set of the neighbors of agent $i$. i.e., $\Nc(i) = \{j \in V \mid \{i,j\}\in E\}$.
The goal of the agents is to cooperatively solve the following optimization problem:
\begin{align}\label{eqn:prob}
\min {}& f(x) \triangleq \sum_{i=1}^m f_i(x)\qquad \text{s.t. } x \in \mathcal{X} \triangleq \bigcap_{i=1}^m \mathcal{X}_i,
\end{align}
where $f_i: \mathbb{R}^d \rightarrow \mathbb{R}$ is a convex function,
representing the local objective of agent $i$,
and $\mathcal{X}_i \subseteq \mathbb{R}^d$  is a closed convex set, representing the local constraint set of agent $i$.
The function $f_i$ and the set $\Xc_i$ are known to agent $i$ only.

We assume that problem~\eqref{eqn:prob} is feasible. Moreover,
we assume each set $\Xc_i$ is defined as the intersection of a collection of simple convex sets.
That is, $\Xc_i$ can be represented as $\Xc_i = \bigcap_{j\in I_i}\Xc_i^j$, where the superscript $j$ is used to identify a component set and $I_i$ is a (possibly infinite) set of indices.
In some applications, $\Xc_i$ may not be explicitly given in advance due to online constraints or uncertainty.
For example, consider the case when $\Xc_i$ is given by
\[
\Xc_i = \{x\in \mathbb{R}^d \mid \la a+\xi ,x\ra \le b\},
\]
where $a\in\mathbb{R}^d$, $b\in\mathbb{R}$ are deterministic and $\xi \in \mathbb{R}^d$ is a Gaussian random noise.
In such a case, a projection-based distributed algorithm cannot be directly applied to solve
problem~(\ref{eqn:prob}) since $|I_i|$ is infinite and the projection of a point on the uncertain set $\Xc_i$ is impossible.
However, a component $\Xc_i^j$ can be realized from a random selection of $\xi$
and the projection onto the realized component is always possible.
Our algorithm is based on such random projections.

We propose a distributed optimization algorithm for problem (\ref{eqn:prob})
that is based on the random projections and the gossip communication protocol.
Gossip algorithms robustly achieve consensus through sparse communications in the network.
That is, only one edge $\{i,j\}$ in the network is randomly selected for communication at each iteration,
and agents $i$ and $j$ simply average their values.
From now on, we refer to our algorithm as \textit{Gossip-based Random Projection} (GRP).

GRP uses an asynchronous time model as in~\cite{Boyd:2006}.
Each agent has a local clock that ticks at a Poisson rate of 1.
The setting can be visualized as having a single virtual clock that ticks whenever any of the local Poisson clock ticks.
Thus, the ticks of the virtual clock is a Poisson random process with rate $m$.
Let $Z_k$ be the absolute time of the $k$th tick of the virtual clock.
The time is discretized according to the intervals $[Z_{k-1},Z_k)$ and this time slot corresponds to our discrete time $k$.
Let $I_k$ denote the index of the agent that wakes up at time $k$ and
$J_k$ denote the index of a neighbor of agent $I_k$ that is selected for communication.
We assume that only one agent wakes up at a time.
The distribution by which $J_k$ is selected is characterized by a nonnegative stochastic $m \times m$ matrix $[\Pi]_{ij} = \pi_{ij}$ that conforms with the graph topology $G =(V,E)$, i.e., $\pi_{ij}>0$ only if
$\{i,j\} \in E$.
At iteration $k$, agent $I_k$ wakes up and contacts one of its neighbors $J_k$ with probability
$\pi_{I_kJ_k}$.

Let $x_{i}(k)$ denote the estimate of agent $i$ at time $k$.
GRP updates these estimates according to the following rule.
Each agent starts with some initial vector $x_i(0)$,
which can be randomly selected.
For $k \ge 1$, agents other than $I_k$ and $J_k$ do not update:
\begin{equation} \label{eqn:algo2}
x_{i}(k) = x_{i}(k-1) \quad \text{for all } i \not\in \{I_k,J_k\}.
\end{equation}
Agents $I_k$ and $J_k$ calculate the average of their estimates, and
adjust the average by using their local gradient information and by projecting onto a randomly selected component of their local constraint sets, i.e.,  for $i \in \{I_k, J_k\}$:
\begin{align}\label{eqn:algo1}
v_{i}(k) {}& = (x_{I_k}(k-1)+x_{J_k}(k-1))/2,\nonumber\\
x_{i}(k) {}& = \mathsf{\Pi}_{\mathcal{X}_i^{\Omega_i(k)}}\left[v_{i}(k) - \alpha_i(k)\nabla f_i(v_{i}(k))\right],
\end{align}
where $\a_i(k)$ is a stepsize of agent $i$, and
$\Omega_i(k)$ is a random variable drawn from the set $I_i$.
The key difference between the work in~\cite{Nedic2009,Nedic2010,KS2011} and this paper is the random projection step.
Instead of projecting on the whole constraint set $\Xc_i$, a component set $\Xc_i^{\Omega_i(k)}$ is selected (or revealed by nature)
and the projection is made on that set,
which reduces the required computations per iteration.

For an alternative representation of GRP we define a nonnegative matrix $W(k)$ as follows:
\[
W(k) = I - \frac{1}{2}(e_{I_k}-e_{J_k})(e_{I_k}-e_{J_k})' \quad \text{ for } k\ge 1,
\]
where $I$ is the $m$-dimensional identity matrix, $e_i \in \mathbb{R}^m$ is a vector whose $i$th entry is equal to 1 and all other entries are equal to 0.
Each $W(k)$ is doubly stochastic by construction, implying that $\Es[W(k)]$ is also doubly stochastic.
Using $W(k)$, algorithm (\ref{eqn:algo2})--(\ref{eqn:algo1}) can be equivalently represented as
\begin{subequations}
\begin{align}
v_{i}(k) {}& = \sum_{j=1}^m[W(k)]_{ij} x_{j}(k-1), \label{eqn:algoa}\\
p_i(k) {}& = \mathsf{\Pi}_{\mathcal{X}_i^{\Omega_i(k)}}[v_{i}(k)-\alpha_i(k)\nabla f(v_{i}(k))]-v_{i}(k),\label{eqn:algob}\\
x_{i}(k) {}& = v_{i}(k) + p_{i}(k)\chi_{\{i\in\{I_k,J_k\}\}},\label{eqn:algoc}
\end{align}
\end{subequations}
where $\chi_{\mathscr{E}}$ is the characteristic-event function, i.e., 
$\chi_{\mathscr{E}}=1$ if $\mathscr{E}$ happens, and $\chi_{\mathscr{E}}=0$ otherwise.

From here onward, we will shorten $\Es[W(k)] = \bar{W}$ since the matrices $W(k)$ are identically distributed.
Let $\lambda$ denote to the second largest eigenvalue of $\bar{W}$.
If  the underlying communication network is connected, 
the incidence graph associated with the positive entries in the matrix
$\bar{W}$ is also connected, with a self-loop at each node. Hence, we have $\lambda <1$.

In the convergence analysis of the algorithm (\ref{eqn:algoa})-(\ref{eqn:algoc}),
we use two different choices of stepsize.
For a diminishing stepsize, we use $\a_i(k) = \frac{1}{\Gamma_i(k)}$ where $\Gamma_i(k)$ denotes the number of updates that agent $i$ has performed until time $k$.
Since every agent $i$ has access to a locally defined quantity $\Gamma_i(k)$,
the stepsize of agent $i$ is independent of every other agent and no coordination is needed for its update.
Another choice that we consider is a constant deterministic stepsize $\a_i(k) = \a_i >0$.

We next discuss our assumptions, the first of which deals with the network.
\begin{assumption}\label{assume:1}
The underlying graph $G =(V,E)$ is connected.
Furthermore, the neighbor selection process is \textit{iid}, whereby at any time agent $i$ is chosen by its neighbor $j \in \Nc(i)$ with probability $\pi_{ji}>0$ ($\pi_{ji}=0$ if $j \not\in \Nc(i)$) independently of the other agents in the network.
\end{assumption}

We use the following assumption for the functions $f_i$ and the sets $\Xc_i^j$.
\begin{assumption} \label{assume:f} Let the following conditions hold:
\begin{enumerate}
\item[(a)] The sets $\mathcal{X}_i^j$, $j\in I_i$, are closed and convex for every $i \in V$.
\item[(b)] Each function $f_i:\mathbb{R}^d\to\mathbb{R}$ is convex over $\mathbb{R}^d$.
\item[(c)] Each function $f_i$ is differentiable and has
\textit{Lipschitz gradients} with a constant $L_i$ over $\mathbb{R}^d$,
\vspace{-2mm}
\[
\|\nabla f_i(x) - \nabla f_i(y)\| \leq L_i \|x-y\| \quad \text{for all } x, y \in \mathbb{R}^d.
\]
\item[(d)] The gradients $\nabla f_i(x)$ are bounded over the set $\mathcal{X}$, i.e.,
there is a constant $G_f$ such that
\vspace{-2mm}
\[
\|\nabla f_i(x)\|  \leq G_f \quad \text{for all $x \in\mathcal{X}$ and all $i\in V$}.
\]
\end{enumerate}
\end{assumption}
For example, Assumption~\ref{assume:f}(d) is satisfied when the constraint set $\mathcal{X}$ is compact.

The next assumption states set regularity, which is crucial in our convergence analysis.
\begin{assumption}\label{assume:c}
There exists a constant $c >0$ such that for all $i\in V$ and $x\in \mathbb{R}^d$,
\begin{equation*}
\mathrm{dist}^2(x,\mathcal{X})
\leq c\, \mathsf{E}\left[\mathrm{dist}^2(x,\mathcal{X}_i^{\Omega_i(k)})|
\Omega_\ell(t),t\in[1,k), \ell\in V\right].
\end{equation*}
Assumption \ref{assume:c} holds if each set $\Xc_i^{j}$ is affine, or the constraint set $\Xc$ has a nonempty interior.
\end{assumption}

\section{Main Results  \label{sec:main}}
In this section, we state the main results of this paper.
The detailed proofs of these results are given later on in Sections~\ref{sec:conv} and~\ref{sec:err}.
We introduce the following notation regarding the optimal value and optimal
solutions of problem~\eqref{eqn:prob}:
\[
f^* = \min_{x\in\mathcal{X}} f(x), \qquad \mathcal{X}^* = \{x \in \mathcal{X} \mid f(x) = f^*\}.
\]
Our first result shows the convergence of the method with probability~1 for a diminishing stepsize.

\begin{proposition} [Convergence w.p.1] \label{prop:p1}
Let Assumptions \ref{assume:1}-\ref{assume:c} hold.
Assume that problem (\ref{eqn:prob}) has a nonempty optimal set $\mathcal{X}^*$ and
the iterates $\{x_i(k)\}$ are generated by algorithm (\ref{eqn:algoa})-(\ref{eqn:algoc})
with $\a_i(k) = 1/\Gamma_i(k)$.
Then, the sequences $\{x_i(k)\}$, for $i \in V$, converge to some random point $x^{\star}$ in the optimal set
$\Xc^*$ with probability 1, i.e., $
\lim_{k\rightarrow \infty} x_i(k) = x^{\star}$ w.p.1 for all  $i\in V.$
\end{proposition}

Proposition~\ref{prop:p1} states that the agents asymptotically reach an agreement
on a random point in the optimal set $\Xc^*$. To get some insights into the convergence rate,
we consider  a constant stepsize $\a_i(k) = \a_i >0$ for $i \in V$, and establish a limiting error bound assuming that each $f_i$ is strongly convex over the set $\Xc$ with a constant $\sigma_i>0$.
The bound will depend on the probabilities of agent updates, which we formally describe as follows.
Let $E_i(k) = \{i \in \{I_k,J_k\}\}$ be the event that agent $i$ updates at time $k$, and
let $\g_i$ be the probability of the event $E_i(k)$. Then,
$\g_i = \frac{1}{m}+\frac{1}{m}\sum_{j \in \Nc(i)} \pi_{ji}$ for all  $i \in V,$
where $\pi_{ji} > 0$ is the probability that agent $i$ is chosen by its neighbor $j$ to communicate.

For the constant stepsize, we will also use the following assumption.
\begin{assumption}\label{assume:s}
Let the convexity requirement for $f_i$ in Assumption~\ref{assume:f}(b) be replaced
by the requirement that each function $f_i$ is
strongly convex with a constant $\sigma_i>0$ over $\mathbb{R}^d$.
In addition, assume that the stepsizes $\a_i$ are such that for all $i\in V$:
\begin{itemize}
\item[(a)] $0<\a_i\sigma_i-4(2+c)\a_i^2 L_i^2 <1$;
\item[(b)] $0<\g_i\left(\a_i\sigma_i-4(2+c)\a_i^2 L_i^2\right) -\frac{\Delta_{\g\a}}{m}<1$,
where $\Delta_{\g\a} = \max_i \{\g_i\a_i\} -\min_j \{\g_j\a_j\}$.
\end{itemize}
\end{assumption}

We have the following result for the asymptotic error bound.
\begin{proposition} [Error bound] \label{prop:sc}
Let Assumptions \ref{assume:1}-\ref{assume:s} hold.
Then, for the iterate sequences $\{x_i(k)\}$, $i \in V$, generated
by algorithm (\ref{eqn:algoa})-(\ref{eqn:algoc}) with a constant stepsize
$\a_i(k) = \a_i >0$, we have
\begin{align*}
\limsup_{k \rightarrow \infty}\frac{1}{m}\sum_{i=1}^m\Es[\|x_i(k)-x^*\|^2]
\le \frac{1}{q}\,4\bar\g\bar \a^2G_f^2\left(\frac{\sqrt{C}}{1-\sqrt{\lambda} } +2(1+c)\right)
+\frac{1}{q}\,\Delta_{\g\a}G_f^2,
\end{align*}
where $x^*$ is the (unique) solution to problem~\eqref{eqn:prob},
\[q=\min_i\{\g_i\rho_i\}- \Delta_{\g\a}/m,\quad\rho_i=\a_i\sigma_i-8(1+c)\a_i^2L_i^2,
\qquad\hbox{for all $i\in V$},\]
$\bar{\g} = \max_i \g_i$, $\bar{\a} = \max_i \a_i$,
 and
$C= 4\left(\frac{8\bar{\g}(1+ \bar{\a}^2\bar L^2)(1+c)}{\min_j\{\g_j\rho_j\} }+ 1\right)$.
\end{proposition}

Proposition \ref{prop:sc} provides an asymptotic error bound for the average of
the expected distances between the iterates of GRP algorithm and
the optimal solution $x^*$.
The first error term is an error term due to a combined effects of the distributed computations over the network, which is controlled by the spectral gap $1-\sqrt{\lambda}$ of the matrix $\bar{W}$,
and the non-diminishing stepsize (common to gradient descent algorithms).
The last term involves an error term $\Delta_{\g\a}G_f^2$
due to the different values for
$\g_i\a_i$ for different agents. We note that if $\g_i\a_i=\nu$ for some $\nu\in(0,1)$ and for all agents $i$,
then this error would be 0.  The condition $\g_i\a_i=\nu$ will hold
when the graph is regular (all $\g_i$ are the same)
and all agents use the same stepsize $\a_i=\a$. There is another more interesting case
when $\g_i\a_i=\nu$ holds for all $i$, which is as follows: the agents that update more frequently
use a smaller stepsize, while the agents that update les frequently use a larger stepsize,
i.e., if $\g_i>\g_j$ then $\a_i<\a_j$, and vice versa.

\section{Convergence Analysis\label{sec:conv}}
In this section,  we prove Proposition~\ref{prop:p1}.
We start with some basic results from the literature, which will be used
later on. The analysis relies on the nonexpansive projection property
 (see \cite{BNO} for its proof), stating that:
for a closed convex set $\mathcal{X} \subseteq \mathbb{R}^d$,
the projection mapping $\mathsf{\Pi}_{\mathcal{X}}: \mathbb{R}^d \rightarrow \mathcal{X}$ is strictly nonexpansive,
\begin{equation}\label{eq:proj}
\|\mathsf{\Pi}_{\mathcal{X}}[x]-y\|^2 \leq \|x-y\|^2 - \|\mathsf{\Pi}_{\mathcal{X}}[x]-x\|^2
\quad \hbox{for all }x \in \mathbb{R}^d \textrm{ and } \hbox{for all } y \in \mathcal{X},
\end{equation}
and, therefore, it is continuous. As an immediate consequence of the preceding relation, we have
\begin{equation}\label{eq:proj-c}
\|\mathsf{\Pi}_{\mathcal{X}}[x]-\mathsf{\Pi}_{\mathcal{X}}[v]\| \leq \|x-v\|
\quad \hbox{for all }x,v \in \mathbb{R}^d.\end{equation}

We also make use of the following
convergence result (see \cite[Lemma 11, p. 49-50]{polyak}).
\begin{lemma}\label{thm:super}
Let $\{v_k\}$, $\{u_k\}$, $\{a_k\}$ and $\{b_k\}$ be non-negative random sequences
such that $
\mathsf{E}[v_{k+1}\mid F_k] \leq (1+a_k)v_k - u_k + b_k$
for all $k \geq 0$ w.p.1,
where $F_k=\{ \{v_i,u_i,a_i,b_i\},0\le i\le k\}$.
If $\sum_{k=0}^{\infty} a_k < \infty$ and $\sum_{k=0}^{\infty} b_k <\infty$ w.p.1,
then  $\lim_{k \rightarrow \infty} v_k = v$ for a random variable $v \geq 0$ w.p.1,
and $\sum_{k=0}^{\infty} u_k < \infty$ w.p.1.
\end{lemma}

The GRP algorithm has three random elements: random gossip communications,
random stepsizes and random projections, which are all independent. They will be handled as follows.

\noindent{\bf Random Gossip Communications:}
At each iteration of the algorithm,
a gossip communication matrix $W(k)$ is realized independently of the past.
In the analysis, we can work with the expected matrix $\bar{W}$ instead of $W(k)$
due to the following properties of the matrices $W(k)$:
(1) Each $W(k)$ is a symmetric projection matrix; hence
$\bar{W}'\bar{W} = \bar{W}$ and
$\left(\bar{W}-\frac{1}{m}\1b\1b'\right)^2 = \bar{W}-\frac{1}{m}\1b\1b'$.
(2) Since $\bar{W}$ is doubly stochastic, the largest eigenvalue of $\bar{W}$ is 1. Therefore, the largest eigenvalue of the matrix $\bar{W}-\frac{1}{m}\1b\1b'$ is the same as
$\lambda$ (the second largest eigenvalue of $\bar{W}$).
These two properties immediately yield the following relation for any
$y \in \mathbb{R}^m$,
\begin{align}\label{eq:ly}
\left\|\left(\bar{W}-\frac{1}{m}\1b\1b'\right)y\right\|^2 \le \lambda\|y\|^2.\end{align}
Furthermore, in view of the connectivity of the underlying graph (Assumption \ref{assume:1}),
we have $\lambda <1$.

\noindent{\bf Random Stepsizes:}
Since the underlying communication graph $G=(V,E)$ is static, due to the gossip-based communications,
the random diminishing stepsize $\a_i(k) = \frac{1}{\Gamma_i(k)}$ exhibits the same behavior as the deterministic stepsize $1/k$ in a long run.
This enables us to handle the cross dependencies of the random stepsizes and the other randomness in the GRP method.

\noindent{\bf Random Projections:}
A projection error is incurred at each iteration of the algorithm since the GRP projects onto
one randomly selected set from the collection defining the overall constraint set $\Xc$.
However, due to the regularity property in the expected sense,
as given in Assumption~\ref{assume:c}, the random projections drive the iterates toward the constraint set $\Xc$ w.p.1 (cf. Lemma~\ref{lemma:key}).

Our convergence analysis is guided by the preceding observations,
and it is constructed along the following main lines:
(1)~the estimates $v_i(k)$ are approaching the constraint set $\Xc$ asymptotically w.p.1;
(2)~the distances $\|v_i(k)-x_i(k)\|$ diminish with probability~1; and
(3)~the agents' estimates $x_i(k)$ eventually arrive at a consensus point that lies in the optimal set $\Xc^*$.
For this, we first establish a basic relation for the iterates of the GRP algorithm
(Lemma~\ref{lem:second}),
which allows us to apply the (almost) supermartingale convergence result of
Lemma~\ref{thm:super}, by letting $v_{k}=\sum_{i=1}^m\|x_i(k)-x^*\|^2$ for some optimal point $x^*$.
To accommodate the use of the  (almost) supermartingale convergence result, we use several auxiliary results.

\subsection{Basic Results for GRP}\label{subsec:iter}
We define the history of the algorithm as follows.
Let $\Fc_k$ be the $\sigma$-algebra generated by the entire history of the algorithm up to time $k$ inclusively, i.e., for all $k \ge 1$,
\[
\Fc_k = \{x_i(0);i\in V\} \cup \{I_{\ell},J_{\ell}, \Omega_i(\ell); i \in \{I_{\ell},J_{\ell}\}, 1 \le \ell \le k\},
\]
and $\Fc_0 = \{x_i(0);i\in V\}$.

We provide several important relations for GRP method.
At first, we provide a relation for the iterates obtained after one step of
algorithm~(\ref{eqn:algoa})-(\ref{eqn:algoc}) and a point in the constraint set $\Xc$.
The lemma relies on the fact that the event that agent $i$ updates at any time is
independent of the past.

\begin{lemma}\label{lem:second} [Basic Iterate Relation]
Let Assumptions~\ref{assume:f}-\ref{assume:c} hold.
Let $\{x_i(k)\}$ be the iterates generated by the algorithm (\ref{eqn:algoa})-(\ref{eqn:algoc}).
Then, for any $q \in (0,1/2)$ there is a sufficiently large $\hat{k}$,
such that with probability 1,
for all $\check{x} \in \Xc$, $k \ge \hat{k}$ and $i \in V$,
\begin{eqnarray*}
&&\EXP{\|x_i(k)-\check{x}\|^2\mid\Fc_{k-1}} \le
 \left(1+\frac{a_1}{k^2}\right)\EXP{\|v_i(k)-\check{x}\|^2 \mid\Fc_{k-1}}
 - \frac{2}{k} \EXP{f_i(z_i(k))-f_i(\check{x})\mid\Fc_{k-1}}\cr
&& -\frac{\g_i}{4c} \EXP{\dist^2(v_i(k),\Xc)\mid\Fc_{k-1}}
+ \frac{a_2}{k^{\frac{3}{2}-q}}
+ \frac{a_3}{k^{\frac{3}{2}-q}} \EXP{\|z_i(k)-\check{x}\|^2\mid\Fc_{k-1}}.
\end{eqnarray*}
where $z_i(k)=\proj_{\Xc}[v_i(k)]$, $a_j>0$ are some constants,
$c$ is the scalar from Assumption~\ref{assume:c},
and $\g_i$ is the probability that agent $i$ updates.
\end{lemma}
The proof of the lemma is in~Appendix \ref{app:lemsecond}, where the constants $a_i$ are also defined.

In the next lemma, we show that the distances between the estimates $v_i(k)$
and the constraint set $\Xc$ go to zero for all $i$, with probability 1 as $k \to \infty$.
We also show that the errors $e_i(k) = x_i(k)-v_i(k)$ converge to zero with probability~1.

\begin{lemma}\label{lemma:key}[Projection Error]
Let Assumptions \ref{assume:f}-\ref{assume:c} hold.
Then, with probability 1, we have
\begin{enumerate}
\item[(a)]
$\sum_{k=1}^\infty \Es\left[\dist^2(v_i(k),\Xc) \mid \Fc_{k-1}\right] <\infty$  and
$\lim_{k\to\infty}\dist(v_i(k),\Xc)=0$ for all $i \in V$.
\item[(b)]
$\sum_{k=1}^\infty \Es[\|e_i(k)\|^2 \mid\Fc_{k-1}]<\infty$ and $\lim_{k\to\infty}\|e_i(k)\|=0$
for all $i \in V$,
where $e_i(k)=x_i(k)-v_i(k)$  for all $i\in V$ and $k\ge1$.
\end{enumerate}
\end{lemma}

Lemma \ref{lemma:key}(a) and Lemma \ref{lemma:key}(b) imply that
$\lim_{k\to\infty} \dist^2(x_i(k),\Xc) = 0$ with probability 1 for all $i \in V$. However, the lemma
does not imply that the sequences $x_i(k)$ converge, nor that their differences $\|x_i(k)-x_j(k)\|$
are vanishing.
A step toward this is provided by the following lemma, which shows
a relation for the agent disagreements on the vectors $v_i(k)$.
\begin{lemma}\label{lem:disagree}
[Disagreement]
Let Assumptions \ref{assume:1}-\ref{assume:f} hold.
Let $\{v_i(k)\}$ be generated by method (\ref{eqn:algoa})-(\ref{eqn:algoc}) with $\a_i(k) = 1/\Gamma_i(k)$ and $\Gamma_i(k)$ being the number of updates that agent $i$
has performed until time~$k$.
Then, for $\bar{v}(k) = \frac{1}{m}\sum_{i=1}^m v_i(k)$ we have
$\sum_{k=1}^{\infty}\frac{1}{k} \Es[\|v_i(k)-\bar{v}(k)\|\mid\Fc_{k-1}]  < \infty$
with probability 1 for all $i \in V$.
\end{lemma}
The proofs of Lemma~\ref{lemma:key} and Lemma~\ref{lem:disagree}
are, respectively, in Appendix \ref{app:lemkey} and~Appendix~\ref{app:lemdisagree}.
\subsection{Proof of Proposition~\ref{prop:p1}}
We assert the convergence of method (\ref{eqn:algoa})-(\ref{eqn:algoc})
using the lemmas established in Section~\ref{subsec:iter}.
Note that Lemma~\ref{lemma:key} allows us to infer that $v_i(k)$ approaches the set $\Xc$,
while Lemma~\ref{lem:disagree} allows us to claim that any two sequences $\{v_i(k)\}$ and $\{v_j(k)\}$
have the same limit points with probability~1.
To claim the convergence of the iterates to an optimal solution, it remains to relate the
limit points of $\{v_i(k)\}$ and the solutions of problem~\eqref{eqn:prob}.
This connection is provided by the iterate relation of Lemma~\ref{lem:second}, supported
by the convergence result in Lemma \ref{thm:super}.
We start the proof by invoking Lemma~\ref{lem:second} stating that
for any $q\in(0,1/2)$, and all $\check{x}\in\Xc$ and $k\ge\hat k$, w.p.1 we have
\begin{eqnarray*}
&&\EXP{\|x_i(k)-\check{x}\|^2\mid\Fc_{k-1}} \le
 \left(1+\frac{a_1}{k^2}\right)\EXP{\|v_i(k)-\check{x}\|^2 \mid\Fc_{k-1}}
 - \frac{2}{k} \EXP{f_i(z_i(k))-f_i(\check{x})\mid\Fc_{k-1}}\cr
&& -\frac{\g_i}{4c} \EXP{\dist^2(v_i(k),\Xc)\mid\Fc_{k-1}}
+  \frac{a_2}{k^{\frac{3}{2}-q} }
+\frac{a_3}{k^{\frac{3}{2}-q}} \EXP{\|z_i(k)-\check{x}\|^2\mid\Fc_{k-1}},
\end{eqnarray*}
where $z_i(k)=\Pi_{\Xc}[v_i(k)]$.
Since $\|z_i(k)-\check{x}\|\le \|v_i(k)-\check{x}\|$
by the non-expansive projection property in~Eq.~\eqref{eq:proj-c},
we obtain
\begin{align}\label{eq:middle}
&\EXP{\|x_i(k)-\check{x}\|^2\mid\Fc_{k-1}} \le
\left(1+\frac{a_4}{k^{\frac{3}{2}-q}}\right) \EXP{\|v_i(k)-\check{x}\|^2 \mid\Fc_{k-1}}\cr
& - \frac{2}{k} \EXP{f_i(z_i(k))-f_i(\check{x})\mid\Fc_{k-1}}
 -\frac{\g_i}{4c}  \EXP{\dist^2(v_i(k),\Xc)\mid\Fc_{k-1}}+\frac{a_2}{k^{\frac{3}{2}-q} },
\end{align}
where $a_4 = a_1+a_3$.
Further, by the definition of $v_i(k) $ in~(\ref{eqn:algoa}), the convexity of the squared-norm function and
the doubly stochastic matrices $W(k)$, we have
\begin{equation}\label{eqn:ds1}
\sum_{i=1}^m \Es[\|v_i(k)-x^*\|^2 \mid \Fc_{k-1}] \le \sum_{i=1}^m\sum_{j=1}^m \bar{W}_{ij}\|x_j(k-1)-x^*\|^2 = \sum_{j=1}^m \|x_j(k-1)-x^*\|^2.\end{equation}
Summing relations in~\eqref{eq:middle} over $i$ and using Eq.~\eqref{eqn:ds1},
yields w.p.1 for all $\check{x}\in\Xc$ and all $k\ge\hat k$,
\begin{eqnarray}\label{eq:propa}
&&\sum_{i=1}^m\EXP{\|x_i(k)-\check{x}\|^2\mid\Fc_{k-1}} \le
 \left(1+\frac{a_4}{k^{\frac{3}{2}-q}}\right)\sum_{j=1}^m\|x_i(k-1)-\check{x}\|^2\cr
 &&- \frac{2}{k} \sum_{i=1}^m\EXP{f_i(z_i(k))-f_i(\check{x})\mid\Fc_{k-1}}
+ \frac{a_2m}{k^{\frac{3}{2}-q} }.
\end{eqnarray}

Recall that $f(x)=\sum_{i=1}^m f_i(x)$. Let $\bar{z}(k) \triangleq \frac{1}{m}\sum_{\ell =1}^m z_\ell(k)$.
Using $\bar{z}(k)$ and $f$,
we can rewrite the term $f_i(z_i(k))-f_i(\check{x})$ as follows:
\begin{align}\label{eqn:rewrite}
\sum_{i=1}^m (f_i(z_i(k))-f_i(\check{x})) {}
& =  \sum_{i=1}^m (f_i(z_i(k))-f_i(\bar{z}(k))) +  (f(\bar{z}(k))-f(\check{x})).
\end{align}
Furthermore, using the convexity of each function $f_i$, we obtain
\begin{align*}
\sum_{i=1}^m (f_i(z_i(k))-f_i(\bar{z}(k)))
 \geq \sum_{i=1}^m  \langle \nabla f_i(\bar{z}(k),z_i(k)-\bar{z}(k)\rangle
\ge -\sum_{i=1}^m  \|\nabla f_i(\bar{z}(k))\|\,\|z_i(k)-\bar{z}(k)\|.
\end{align*}
Since $\bar{z}(k)$ is a convex combination of $z_i(k)\in \Xc$, it follows that
$\bar z(k) \in \Xc$. Using $\bar z(k) \in \Xc$ and the uniform bound $G_f$
for the norms $\|\nabla f_i(x)\|$ on the set $\Xc$ (Assumption~\ref{assume:f}(d))
we obtain
\begin{align}\label{eqn:a}
\sum_{i=1}^m (f_i(z_i(k))-f_i(\bar{z}(k)))
\ge -G_f\sum_{i=1}^m \|z_i(k)-\bar{z}(k)\|.
\end{align}
 We next consider the term $\|z_i(k)-\bar{z}(k)\|$, for which by using
 $\bar{z}(k) \triangleq \frac{1}{m}\sum_{\ell =1}^m z_\ell(k)$ we have
 \[\|z_i(k)-\bar{z}(k)\|=\left\|\frac{1}{m}\sum_{\ell=1}^m(z_i(k)-z_\ell(k))\right\|
 \le \frac{1}{m}\sum_{\ell=1}^m\|z_i(k)-z_\ell(k)\|
 \le \frac{1}{m}\sum_{\ell=1}^m\|v_i(k)-v_\ell(k)\|,
 \]
 where the first inequality is obtained by the convexity of the norm
 and the last inequality follows by the projection property in Eq.~\eqref{eq:proj-c}.
 Further, by letting $\bar v(k)=\frac{1}{m}\sum_{j=1}^m v_j(k)$ and
 using $\|v_i(k)-v_\ell(k)\|\le \|v_i(k)-\bar v(k)\|+ \|v_\ell(k)- \bar v(k)\|$,
 we obtain
 $\|z_i(k)-\bar{z}(k)\|\le \|v_i(k)-\bar v(k)\|+\frac{1}{m}\sum_{\ell=1}^m\|v_\ell(k)-\bar v(k)\|$
 for every $i\in V$.
Upon summing these relations over $i\in V$, we find
\begin{align}\label{eqn:sums}
\sum_{i=1}^m \|z_i(k)-\bar{z}(k)\|\le 2 \sum_{i=1}^m\|v_i(k)-\bar v(k)\|.
\end{align}
Combining relations~\eqref{eqn:sums} and~\eqref{eqn:a},
and substituting the resulting relation in Eq.~\eqref{eqn:rewrite}, we obtain
\begin{align}\label{eqn:z_i}
\sum_{i=1}^m (f_i(z_i(k))-f_i(\check{x}))
\ge -2G_f \sum_{i=1}^m \|v_i(k)-\bar{v}(k)\|+  (f(\bar{z}(k))-f(\check{x})).
\end{align}
Finally, by using the preceding estimate in inequality~\eqref{eq:propa}
and letting $\check{x}=x^*$ for an arbitrary $x^* \in \Xc^*$,
we have w.p.1 for any $x^* \in \mathcal{X}^*$ and $k \geq \bar{k}$,
\begin{eqnarray}\label{eq:propb}
&&\sum_{i=1}^m\EXP{\|x_i(k) -x^*\|^2\mid\Fc_{k-1}} \le
 \left(1+\frac{a_4}{k^{\frac{3}{2}-q}}\right)
 \sum_{j=1}^m\|x_i(k-1) - x^*\|^2\cr
 &&- \frac{2}{k} \EXP{f(\bar z(k))-f^*)\mid\Fc_{k-1}}
 +\frac{4G_f}{k} \sum_{i=1}^m \Es[\|v_i(k)-\bar{v}(k)\|\mid\Fc_{k-1}]
+ \frac{a_2m}{k^{\frac{3}{2}-q} }.\qquad
\end{eqnarray}
Since $\bar{z}(k) \in \Xc$, we have $f(\bar{z}(k)) - f^* \ge 0$.
Thus, in the light of Lemma \ref{lem:disagree}, relation~(\ref{eq:propb}) satisfies all the conditions of
Lemma~\ref{thm:super}.
Hence, the sequence $\{\|x_i(k)-x^*\|^2\}$ is convergent
for any $i \in V$ and $x^* \in \Xc^*$ w.p.1, and
$\sum_{k=0}^{\infty} \frac{1}{k}(f(\bar{z}(k))-f^*) < \infty$ w.p.1.
Since $\sum_{k=0}^{\infty} \frac{1}{k} = \infty$, it follows that
\begin{equation}\label{eqn:res1}
\liminf_{k\rightarrow \infty} (f(\bar{z}(k))-f^*)=0\quad \text{ w.p.1.}
\end{equation}
By Lemma~\ref{lemma:key}(a), noting that $z_i(k)=\mathsf{\Pi}_{\Xc}[v_i(k)]$, we have
\begin{equation}\label{eqn:res2}
\lim_{k\rightarrow \infty}\|v_i(k)-z_i(k)\|=0 \quad \text{for all } i \in V \text{ w.p.1}.
\end{equation}
Since the sequence $\{\|x_i(k)-x^*\|\}$ is convergent with probability 1 for any $i \in V$ and every $x^* \in \Xc^*$, in view of the relations (\ref{eqn:algoa}) and (\ref{eqn:res2}), respectively,
so are the sequences $\{\|v_i(k)-x^*\|\}$ and $\{\|z_i(k)-x^*\|\}$, as well as their
average sequences $\{\|\bar v(k)-x^*\|\}$ and $\{\|\bar{z}(k)-x^*\|\}$.
Therefore, the sequences $\{\bar v(k)\}$ and $\{\bar z(k)\}$
are bounded with probability 1, and they have accumulation points.
From relation (\ref{eqn:res1}) and the continuity of $f$, the sequence $\{\bar{z}(k)\}$ must have one accumulation point in $\Xc^*$ with probability 1.
This and the fact that $\{\|\bar{z}(k)-x^*\|\}$ is convergent with probability 1 for every $x^* \in \Xc^*$ imply that for a random point $x^{\star} \in \Xc^*$,
\begin{equation}\label{eqn:z_final}
\lim_{k \rightarrow \infty} \bar{z}(k) = x^{\star} \quad~ \text{ w.p.1}.
\end{equation}

Now, from $\bar{z}(k) = \frac{1}{m}\sum_{\ell=1}^m z_\ell(k)$
and $\bar{v}(k) = \frac{1}{m}\sum_{i=\ell}^m v_\ell(k)$, using
relation~(\ref{eqn:res2}) and the convexity of the norm, we obtain
$\lim_{k\to\infty}\|\bar{v}(k)-\bar{z}(k)\|\le\frac{1}{m}\sum_{\ell=1}^m
\lim_{k\to\infty}\|v_\ell(k)-z_\ell(k)\|=0$ w.p.1.
In view of relation~\eqref{eqn:z_final}, it follows that
\begin{equation}\label{eqn:v_final}
\lim_{k \rightarrow \infty} \bar{v}(k) = x^{\star} \quad~ \text{ w.p.1}.
\end{equation}
By  Lemma \ref{lem:disagree}, we have
\begin{align}\label{eqn:liminfv}
\liminf_{k \rightarrow \infty} \|v_i(k)-\bar{v}(k)\| = 0 \quad \text{for all } i \in V \text{ w.p.1}.
\end{align}
The fact that $\{\|v_i(k)-x^*\|\}$ is convergent with probability 1 for all $i$ and any $x^* \in \Xc^*$, together with~\eqref{eqn:v_final}
and~\eqref{eqn:liminfv}  implies that
\begin{equation}\label{eqn:consensus}
\lim_{k \rightarrow \infty} \|v_i(k)-x^\star\| = 0 \quad \text{for all } i \in V \text{ w.p.1}.
\end{equation}
Finally, by Lemma \ref{lemma:key}(b), we have
$\lim_{k \rightarrow \infty} \|x_i(k)-v_i(k)\| = 0$ for all $i \in V$ w.p.1,
which together with the limit in~\eqref{eqn:consensus} yields
$\lim_{k\rightarrow \infty} x_i(k) = x^\star$ for all $i \in V$ with probability 1.\quad$\blacksquare$

\section{Error Bound\label{sec:err}}
Here, we prove Proposition~\ref{prop:sc}.
We start by providing some lemmas that are valid for a constant stepsize $\a_i(k) = \a_i >0$.
The first result shows a basic iterate relation.

\begin{lemma}\label{lem:base}
Let Assumptions \ref{assume:f}-\ref{assume:s} hold, where
the stepsize satisfies Assumption~\ref{assume:s}(a).
Then, for the iterates $x_i(k)$ of the method we have w.p.1
for any $x\in \Xc$, and for all $k\ge1$ and $i \in \{I_k,J_k\}$,
\begin{align*}
\EXP{\|x_i(k)-x\|^2 \mid \Fc_{k-1},I_k,J_k}
\leq{}& (1-\rho_i) \|v_i(k)-x\|^2
- 2\alpha_i\langle \nabla f_i(x), z_i(k)-x\rangle + 8(1+c)\a_i^2G_f^2, 
\end{align*}
where $\rho_i=\sigma_i\alpha_i-4(2+c)\a_i^2L_i^2$ and $z_i(k)=\Pi_{\Xc}[v_i(k)]$.
\end{lemma}
The proof of the preceding lemma is in Appendix~\ref{app:lembase}.

Next, we provide an asymptotic estimate for the disagreement among the agents.
\begin{lemma}\label{lem:disagree2}
Let Assumptions \ref{assume:1}--\ref{assume:s} hold.
Let $\bar{x}(k) = \frac{1}{m}\sum_{i=1}^m x_i(k)$ for all $k$.
Then, for the iterates $\{x_i(k)\}$ generated by method (\ref{eqn:algoa})-(\ref{eqn:algoc}), we have
\[
\limsup_{k \rightarrow \infty}\sum_{i=1}^m \Es[\|x_i(k)-\bar{x}(k)\|^2]
\le \frac{4m\bar{\a}^2G_f^2}{(1-\sqrt{\lambda})^2}C,
\]
where
$C=\frac{8\bar{\g}(1+ \bar{\a}^2\bar L^2)(1+c)}{\min_j\{\g_j\rho_j\}}+ 1$,
$\rho_i=\sigma_i\alpha_i-4(2+c)\a_i^2L_i^2$,
$\bar{\g} = \max_i \g_i$, $\bar{\a} = \max_i \a_i$, and $\bar L=\max L_i$.
\end{lemma}
The proof of the lemma is given in Appendix~\ref{app:lemdisagree2}.
The bound in Lemma \ref{lem:disagree2} captures the variance of the estimates $x_i(k)$
in terms of the number of agents, the maximum stepsize and the spectral gap $1-\sqrt{\lambda}$ of the matrix $\bar{W}$.

We are now ready to prove Proposition~\ref{prop:sc}. In the proof we use a
relation implied by the convexity of the squared-norm. In particular,
by the definition of $v_i(k) $ in~(\ref{eqn:algoa}), the convexity of the squared-norm function and
the doubly stochastic weights $W(k)$, we have for any $x\in\mathbb{R}^d$,
\begin{equation}\label{eqn:ds2}
\sum_{i=1}^m \Es[\|v_i(k)-x\|^2 \mid \Fc_{k-1}]
\le \sum_{i=1}^m\sum_{j=1}^m \bar{W}_{ij}\|x_j(k-1)-x\|^2 = \sum_{j=1}^m \|x_j(k-1)-x\|^2.
\end{equation}
{\it Proof of Proposition \ref{prop:sc}.}
The function $f$ is strongly convex with a constant $\sigma = \sum_{i=1}^m \sigma_i$ and therefore,
problem (\ref{eqn:prob}) has a unique optimal solution $x^*$.
The proof starts with the relation of Lemma~\ref{lem:base} where we let $x=x^*$.
Define $\bar{z}(k) = \frac{1}{m} \sum_{i=1}^m z_i(k)$, so that $\bar{z}(k) \in \Xc$.
We have
$\langle\nabla f_i(x^*),z_i(k)-x^*\rangle
= \langle\nabla f_i(x^*),\bar{z}(k)-x^*\rangle +
\langle\nabla f_i(x^*),z_i(k)-\bar{z}(k)\rangle$, which in view of
the gradient boundedness (Assumption~\ref{assume:f}(d)) implies that
\[
\langle\nabla f_i(x^*),z_i(k)-x^*\rangle \ge
\langle\nabla f_i(x^*),\bar{z}(k)-x^*\rangle - G_f\|z_i(k)-\bar{z}(k)\|.
\]
Using the preceding relation and Lemma~\ref{lem:base},
we have for all $k \ge 1$ w.p.1,
\begin{align*}
\Es[\|x_i(k)-x^*\|^2\mid \Fc_{k-1},I_k,J_k]
\leq {}& (1-\rho_i)\|v_i(k)-x^*\|^2 -2\alpha_i\langle \nabla f_i(x^*),\bar{z}(k)-x^*\rangle\cr
{}& +2\alpha_iG_f\|z_i(k)-\bar{z}(k)\| + 8(1+c)\alpha_i^2G_f^2.
\end{align*}
Taking the expectation with respect to $\Fc_{k-1}$ and using the fact that the preceding inequality holds with probability $\g_i$, and otherwise we have $x_i(k) = v_i(k)$ with probability $1-\g_i$, we obtain w.p.1 for all
$k \ge 1$ and $i \in V$,
\begin{align*}
&\Es[\|x_i(k)-x^*\|^2\mid \Fc_{k-1}]
\leq  (1-\g_i\rho_i)\Es[\|v_i(k)-x^*\|^2\mid \Fc_{k-1}] \\
{}& -2\g_i\alpha_i\Es[\langle\nabla f_i(x^*),\bar{z}(k)-x^*\rangle\mid \Fc_{k-1}]
+2\g_i\alpha_iG_f\Es[\|z_i(k)-\bar{z}(k)\|\mid \Fc_{k-1}] + 8(1+c)\g_i\alpha_i^2G_f^2.
\end{align*}
We note that under the assumption that
$\rho_i=\a_i\sigma_i-8(1+c)\a_i^2L_i^2\in(0,1)$ for all $i$,  we also have
$\g_i\rho_i\in(0,1)$ for all $i$ since $\g_i\in(0,1)$.
By adding and subtracting
$2\min_j\{\g_j\a_j\}\Es[\langle\nabla f_i(x^*),\bar{z}(k)-x^*\rangle\mid \Fc_{k-1}]$,
we find that
\begin{align}\label{eqn:mid}
&\Es[\|x_i(k)-x^*\|^2\mid \Fc_{k-1}]
\leq  (1-\g_i\rho_i)\,\Es[\|v_i(k)-x^*\|^2\mid \Fc_{k-1}] \cr
{}& -2\underline{\g}\underline{\a}\Es[\langle\nabla f_i(x^*),\bar{z}(k)-x^*\rangle\mid \Fc_{k-1}]
+2\Delta_{\g\a}\Es[\|\nabla f_i(x^*)\| \|\bar{z}(k)-x^*\|\mid\Fc_{k-1}]\cr
{}&+2\bar\g\bar\alpha G_f\Es[\|z_i(k)-\bar{z}(k)\|\mid \Fc_{k-1}] + 8(1+c)\bar\g\bar\alpha^2G_f^2,
\end{align}
where $\Delta_{\g\a} = \max_j\{\g_j\a_j\}-\min_j\{\g_j\a_j\}$,
$\underline{\a} =\min_i \a_i$, $\underline{\g} = \min_i g_i$, $\bar{\a} =\max_i \a_i$ and
$\bar{\g} =\max_i \g_i$.
We can further estimate
\begin{align*}
\|\nabla f_i(x^*)\| \|\bar{z}(k)-x^*\| \le  \frac{G_f}{m} \sum_{i=1}^m\|\mathsf{\Pi}_{\Xc}[v_i(k)]-x^*\|  \le \frac{G_f}{m} \sum_{i=1}^m\|v_i(k)-x^*\|,
\end{align*}
where the first inequality follows by Assumption \ref{assume:f}(d),
$\bar{z}(k) = \frac{1}{m}\sum_{i=1}^m \mathsf{\Pi}_{\Xc} [v_i(k)]$ and the convexity of the norm function,
while the second inequality follows from the projection property~(\ref{eq:proj-c}).
Also, from relation $ab \le \frac{1}{2}(a^2+b^2)$ and the convexity of the square-function,
we obtain
\begin{align}\label{eqn:mid2}
\|\nabla f_i(x^*)\| \|\bar{z}(k)-x^*\|
\le \frac{1}{2}G_f^2+ \frac{1}{2m} \sum_{i=1}^m\|v_i(k)-x^*\|^2.
\end{align}

Summing relations in Eq.~(\ref{eqn:mid}) over $i$, and using estimates (\ref{eqn:ds2}), (\ref{eqn:mid2}) and
$\sum_{i=1}^m \langle\nabla f_i(x^*),\bar{z}(k)-x^*\rangle \ge f(\bar z(k))-f(x^*)\ge 0$
(which holds by the optimality of $x^*$),
we have
\begin{align*}
\sum_{i=1}^m\Es[\|x_i(k)-x^*\|^2\mid \Fc_{k-1}]
\leq {}& (1-q)\sum_{j=1}^m\Es[\|x_j(k-1)-x^*\|^2\mid \Fc_{k-1}]
+\Delta_{\g\a}mG_f^2\cr
{}& +2\bar{\g}\bar{\alpha}G_f\sum_{i=1}^m\Es[\|z_i(k)-\bar{z}(k)\|\mid \Fc_{k-1}]
+ 8(1+c)m\bar{\g}\bar{\alpha}^2G_f^2.
\end{align*}
where $q =\min_i\{\g_i\rho_i\}- \Delta_{\g\a}/m$.
Since $\g_i\rho_i- \frac{\Delta_{\g\a}}{m}\in(0,1)$ by Assumption~\ref{assume:s}(b), it follows that
$q \in(0,1)$, and therefore
\begin{align}\label{eq:last0}
\limsup_{k \rightarrow \infty}\sum_{i=1}^m\Es[\|x_i(k)-x^*\|^2]
\leq{}&
2\bar{\g}\bar{\alpha}G_f\frac{1}{q}\limsup_{k \rightarrow \infty}\sum_{i=1}^m\Es[\|z_i(k)-\bar{z}(k)\|] \cr
{}&+\frac{1}{q}\left(\Delta_{\g\a}mG_f^2+ 8(1+c)m\bar{\g}\bar{\alpha}^2G_f^2\right).
\end{align}

We now consider the sum $\sum_{i=1}^m\Es[\|z_i(k)-\bar{z}(k)\|] $.
Using H\"older's inequality, we have
\begin{align}\label{eq:last1}
\sum_{i=1}^m\Es[\|z_i(k)-\bar{z}(k)\|] \le \sqrt{m\sum_{i=1}^m \EXP{\|z_i(k)-\bar{z}(k)\|^2}}.\end{align}
Since $\bar{z}(k) = \frac{1}{m}\sum_{i=1}^m \mathsf{\Pi}_{\Xc} [v_i(k)]$, it follows that
for $\bar v(k)=\frac{1}{m}\sum_{j=1}^m v_j(k)$,
\begin{align}\label{eq:last2}
\sum_{i=1}^m\Es[\|z_i(k)-\bar{z}(k)\|^2] \le
\sum_{i=1}^m\Es[\|z_i(k)-\mathsf{\Pi}_{\Xc}[\bar{v}(k)]\|^2]
\le \sum_{i=1}^m\Es[\|v_i(k)-\bar{v}(k)\|^2],
\end{align}
where the last inequality is obtained from the projection property~\eqref{eq:proj-c}.
Since $\bar v(k)$ is the average of $v_j(k)$ for $j\in V$, it follows that for
$\bar{x}(k-1)=\frac{1}{m}\sum_{j=1}^m x_j(k-1)$,
\[\sum_{i=1}^m\EXP{\|v_i(k)-\bar{v}(k)\|^2} \le \sum_{i=1}^m \EXP{\|v_i(k)-\bar{x}(k-1)\|^2}.\]
From the preceding relation, and using Eq.~\eqref{eqn:ds2} with $x=\bar{x}(k-1)$ (where we take the
total expectation), we find that
\begin{align}\label{eq:last3}
\sum_{i=1}^m\EXP{\|v_i(k)-\bar{v}(k)\|^2}
\le \sum_{j=1}^m \EXP{\|x_j(k-1)-\bar{x}(k-1)\|^2}.
\end{align}
From Eqs.~\eqref{eq:last1}--\eqref{eq:last3} and Lemma~\ref{lem:disagree2}, we obtain
\[\limsup_{k \rightarrow \infty}\sum_{i=1}^m\Es[\|z_i(k)-\bar z(k)\|]
\le 2 m\sqrt{C}\,\frac{\bar\a G_f}{1-\sqrt{\lambda} }.\]
The result follows from Eq.~\eqref{eq:last0} after dividing by $m$.
\quad$\blacksquare$

\section{Simulations: Distributed Robust Control \label{sec:sim}}
In this section, we apply our GRP algorithm to a distributed robust model predictive control (MPC) problem \cite{Camponogara}.
A linear, time-invariant, discrete-time system is given by the following state equation for $t = 1, \ldots, T,$
\begin{equation}\label{eqn:system}
x(t) = Ax(t-1) + Bu(t),
\end{equation}
where
\[
A =\left[\begin{array}{cc}
1 & 1\\
0 & 1
\end{array}\right],\quad
b = \left[\begin{array}{c}
0.5\\
1
\end{array}\right],
\]
with initial state $x(0) = [7,~0]'$.
The goal of the agents on the network is to find an optimal control
$\mathbf{u} \triangleq [u(1),\ldots,u(T)]'$ of system (\ref{eqn:system}) over time $t = 1,\ldots,T,$ with some random terminal constraints.
The distributed optimization problem is given by
\begin{equation}\label{eqn:mpc}
\min_{\mathbf{u}}~  f(\mathbf{u}) =  \sum_{i=1}^m f_i(\mathbf{u}) \quad \text{s.t. } \mathbf{u} \in \Xc,
\end{equation}
where
\[
f_i(\mathbf{u}) = \sum_{t=1}^T \|x(t)-z_i\|^2 + r u(t), \text{ for } i = 1, \ldots, m,
\]
is the local objective of agent $i$ and $r >0$ is a control parameter.
Hence, the agents on the network jointly find a control $\mathbf{u}$ such that the resulting trajectory
$x(t)$, for $t = 1, \ldots, T,$ minimizes the deviations from the points
$z_i \in \mathbb{R}^2$ together with the control effort.
The information about the points $z_i$, for $i = 1,\ldots,m,$ are private and only agent $i$ knows
the location $z_i$.

The constraint set $\Xc$ is a set of control inputs that satisfies the following constraints.
\begin{subequations}
\begin{align}
&\|u(t)\|_{\infty} \le 2, \quad \text{for } t = 1, \ldots, T, \label{eqn:cons1}\\
&x(t) = Ax(t-1) + Bu(t), \quad \text{for } t = 1, \ldots, T,\label{eqn:cons2}\\
&\max_{\ell=1,2,3,4}\left\{(a_{\ell}+\delta_{\ell})'x(T) -b_{\ell}\right\} \le 0. \label{eqn:cons4}
\end{align}
\end{subequations}
The system is initiated in state $x(0) = [7,~0]'$.
The constraint (\ref{eqn:cons1}) is just a box constraint, while
the constraints in~(\ref{eqn:cons2}) describe the system dynamics.
The constraints in
(\ref{eqn:cons4}) describe the random terminal conditions given by the linear inequalities $(a_{\ell}+\delta_{\ell})'x(T)\le b_{\ell}$
and the perturbations $\delta_{\ell}$ are uniform random vectors in boxes $\|\delta_{\ell}\|_{\infty} \le \b_{\ell}$ for some given scalars $\beta_\ell$.
Note that $u(t)$, for $t = 1, \ldots, T$, are the only variables here since $x(t)$, for $t = 1, \ldots, T,$ are fully determined by state equations~(\ref{eqn:cons2}) once $u(t)$, for $t = 1, \ldots, T,$ is given.

For this problem, we have $\Xc_i = \Xc$ for all $i$.
The constraint set $\Xc$ is uncertain and not exactly known in advance since the perturbations are uniform random vectors in boxes.
To apply the GRP algorithm (\ref{eqn:algoa})-(\ref{eqn:algoc}) in solving this robust optimal control problem, at iteration $k$, each agent $I_k$ and $J_k$  draws a realization of one of the linear inequality terminal constraints,
and each of them projects its current iterate on the selected constraint.
Subsequently, they perform their projections onto the box constraint~(\ref{eqn:cons1}).

Since the uncertainty exists in a box, the problem (\ref{eqn:mpc}) has an equivalent Quadratic Programming (QP) formulation.
Note that the following representations are all equivalent:
\begin{subequations}
\begin{equation}\label{eqn:suba}
(a_{\ell}+\delta_{\ell})'x(T)\le b_{\ell},\quad \forall (\delta_{\ell}: \|\delta_{\ell}\|_{\infty} \le \b_{\ell})
\end{equation}
\begin{equation}\label{eqn:subb}
\Leftrightarrow \max_{\|\delta_{\ell}\|_{\infty} \le \b_{\ell}} \delta_{\ell}'x(T) \le b_{\ell} - a_{\ell}'x(T)
\end{equation}
\begin{equation}\label{eqn:subc}
\Leftrightarrow a_{\ell}'x(T) + \b_{\ell}|[x(T)]_1| + \b_{\ell}|[x(T)]_2|\le b_{\ell}.
\end{equation}
\end{subequations}
Therefore, the inequality (\ref{eqn:cons4}) admits an equivalent representation
of (\ref{eqn:subc}) by a system of linear inequalities with additional variables $t_1$ and $t_2$:
\begin{subequations}
\begin{equation}\label{eqn:qpa}
-t_j \le [x(T)]_j \le t_j,\quad \text{for } j = 1,2,
\end{equation}
\begin{equation}\label{eqn:qpb}
\max_{\ell = 1,2,3,4}\left\{a_{\ell}'x(T) + \b_{\ell}t_1 + \b_{\ell}t_2 - b_{\ell}\right\} \le 0.
\end{equation}
\end{subequations}
This alternative representation is only available since we are considering simple box uncertainty sets
for the sake of comparison. Note that our GRP algorithm is applicable not just to box uncertainty
but to more complicated perturbations such as Gaussian or other distributions.

\begin{figure}[t]
\begin{center}
\includegraphics[scale=0.4]{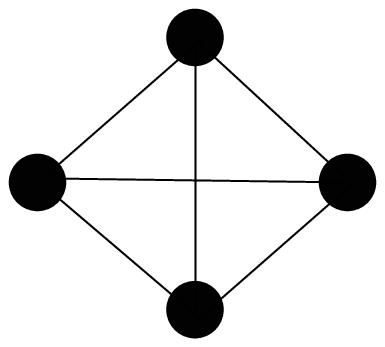}\hspace{0.25in}
\includegraphics[scale=0.4]{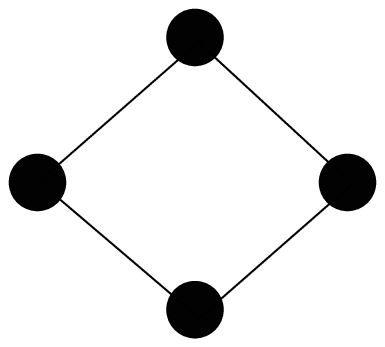}\hspace{0.25in}
\includegraphics[scale=0.4]{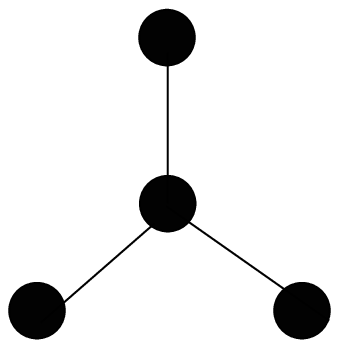}
\end{center}
\caption{\label{fig:top} Clique (left), cycle (center) and star (right) graph used for communication topology ($m=4$)}
\end{figure}

In the experiment, we use $m=4$ and $m=10$ agents with $T=10$ and $r=0.1$.
We solve the problem on three different network topologies, namely, clique, cycle and star
(see Figure~\ref{fig:top}).
For the agent selection probability, we use uniform distribution, i.e.,
at each iteration, one of the $m$ agents is uniformly selected and the selected agent uniformly selects one of its neighbors.
Table~\ref{tbl:lambda} shows the second largest eigenvalue $\lambda$ of $\bar{W}$ for the three network topologies when $m=4$ and $m=10$.
When $m$ is larger, we can see that $\lambda$ is very close to one for all of the three cases.

We evaluate the algorithm performance by carrying out 100 Monte-Carlo runs, each with 40,000 iterations for $m=4$ and 100,000 iterations for $m=10$.
For the stepsize, we use either a diminishing one ($1/\Gamma_i(k)$) or a constant $\alpha_i = 10^{-5}$ for $m=4$ and $\alpha_i = 10^{-6}$ for $m=10$.

\begin{figure}[t]
\begin{center}
\includegraphics[scale=0.35]{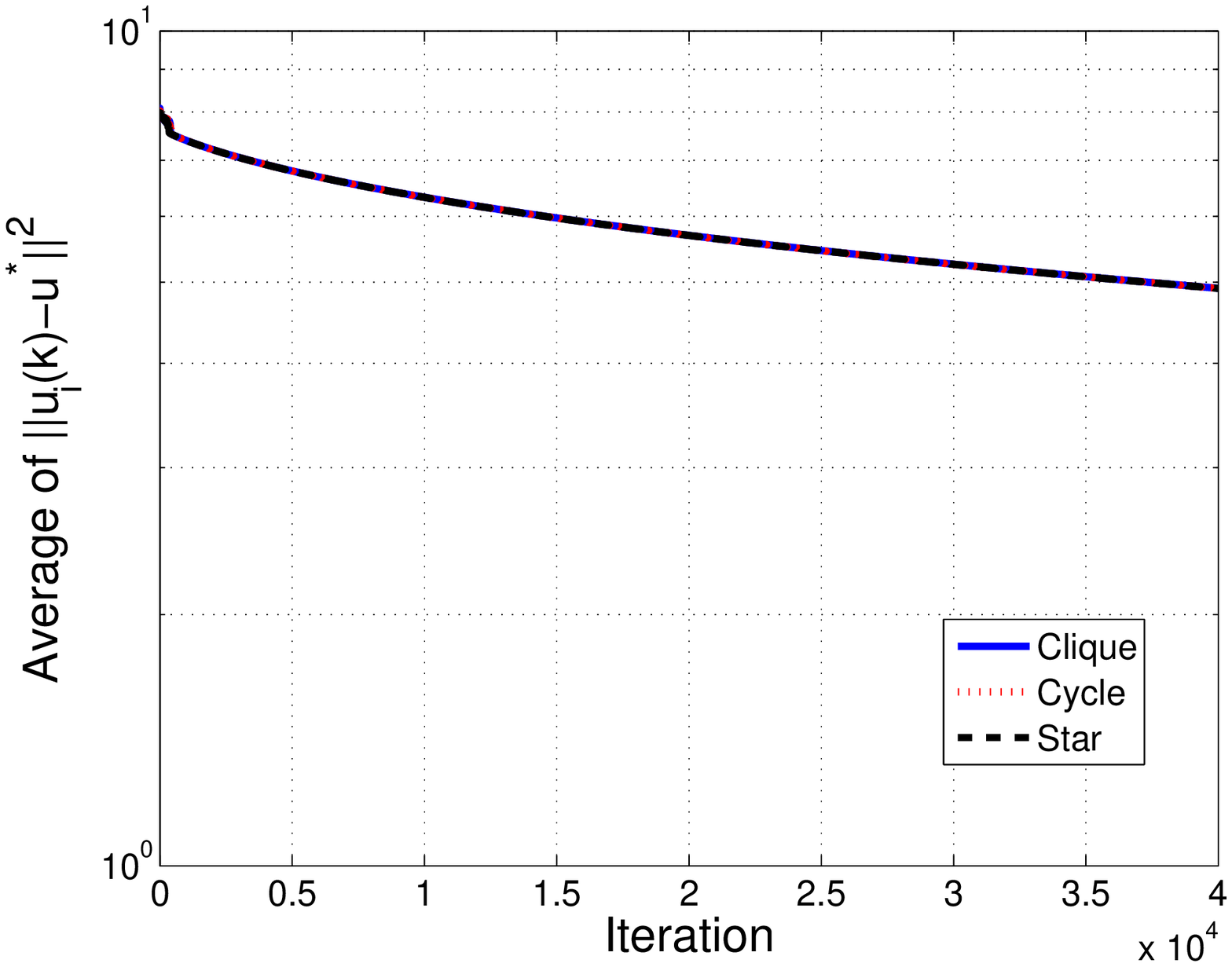} \hspace{0.1in}
\includegraphics[scale=0.35]{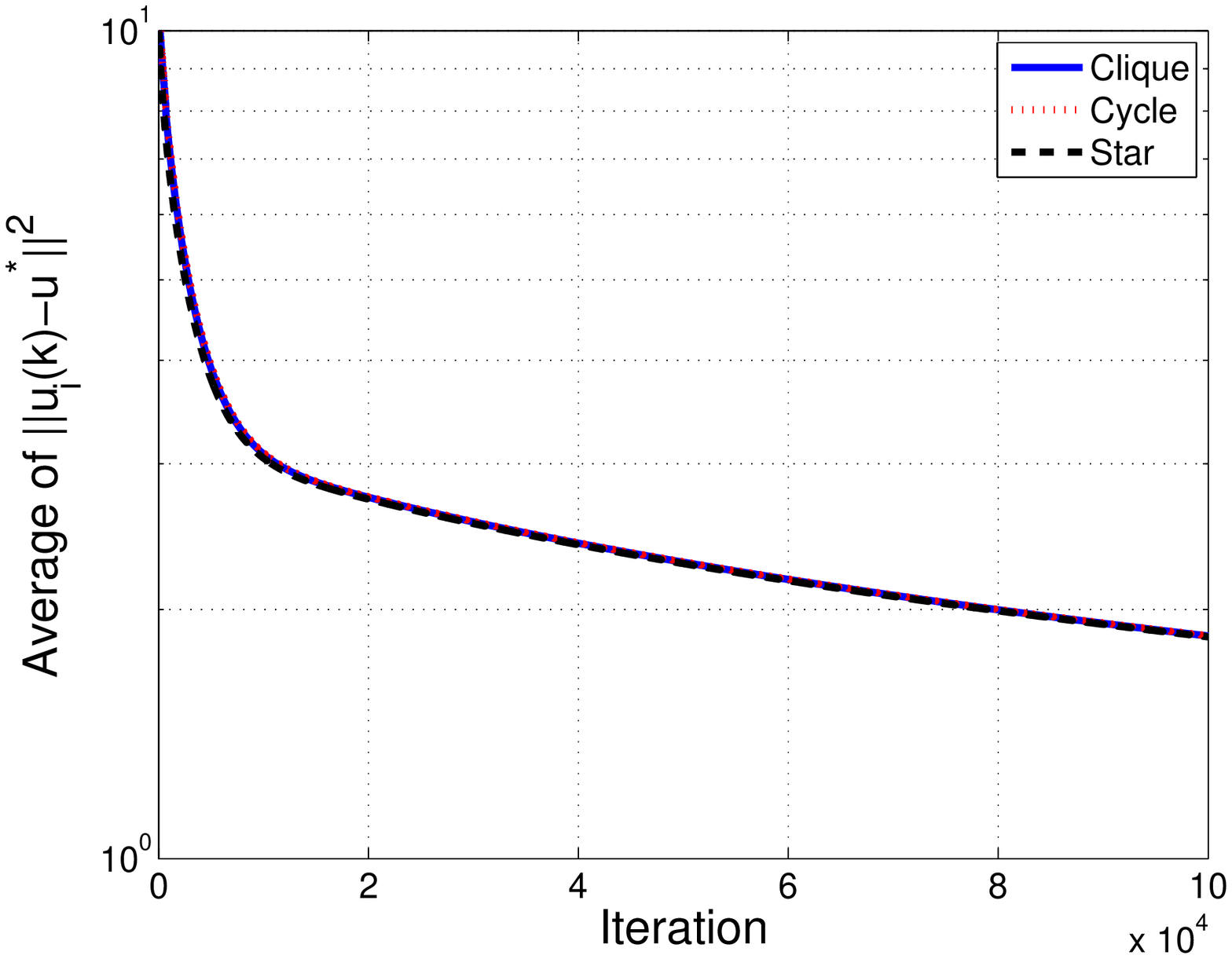}
\end{center}
\caption{\label{fig:m4-dim} Iteration vs $\frac{1}{m}\sum_{i=1}^m\|\mathbf{u}_i(k)-\mathbf{u}^*\|^2$ with a diminishing stepsize when $m=4$ (left) and $m=10$ (right)}
\end{figure}

\begin{figure}[t]
\begin{center}
\includegraphics[scale=0.35]{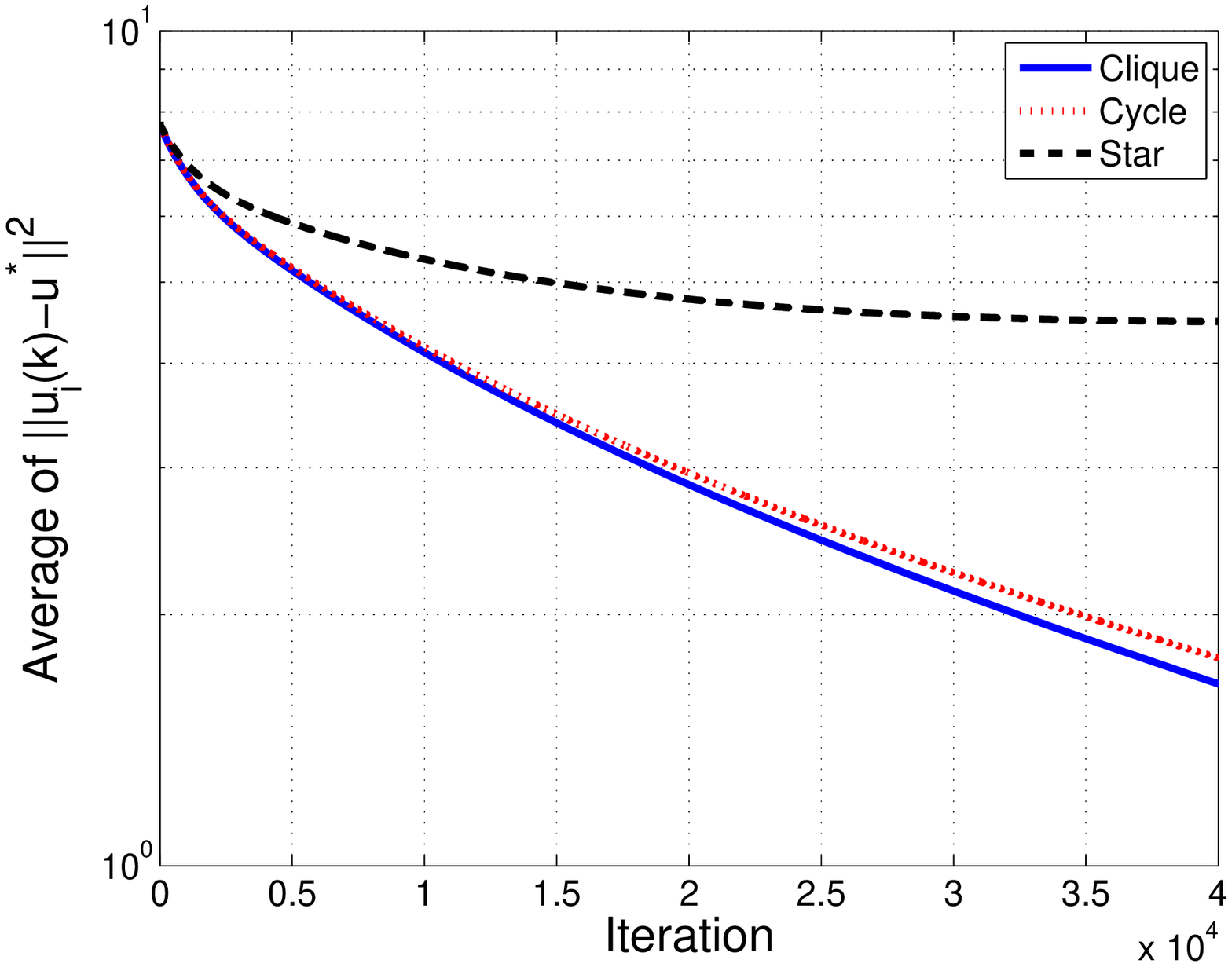} \hspace{0.1in}
\includegraphics[scale=0.35]{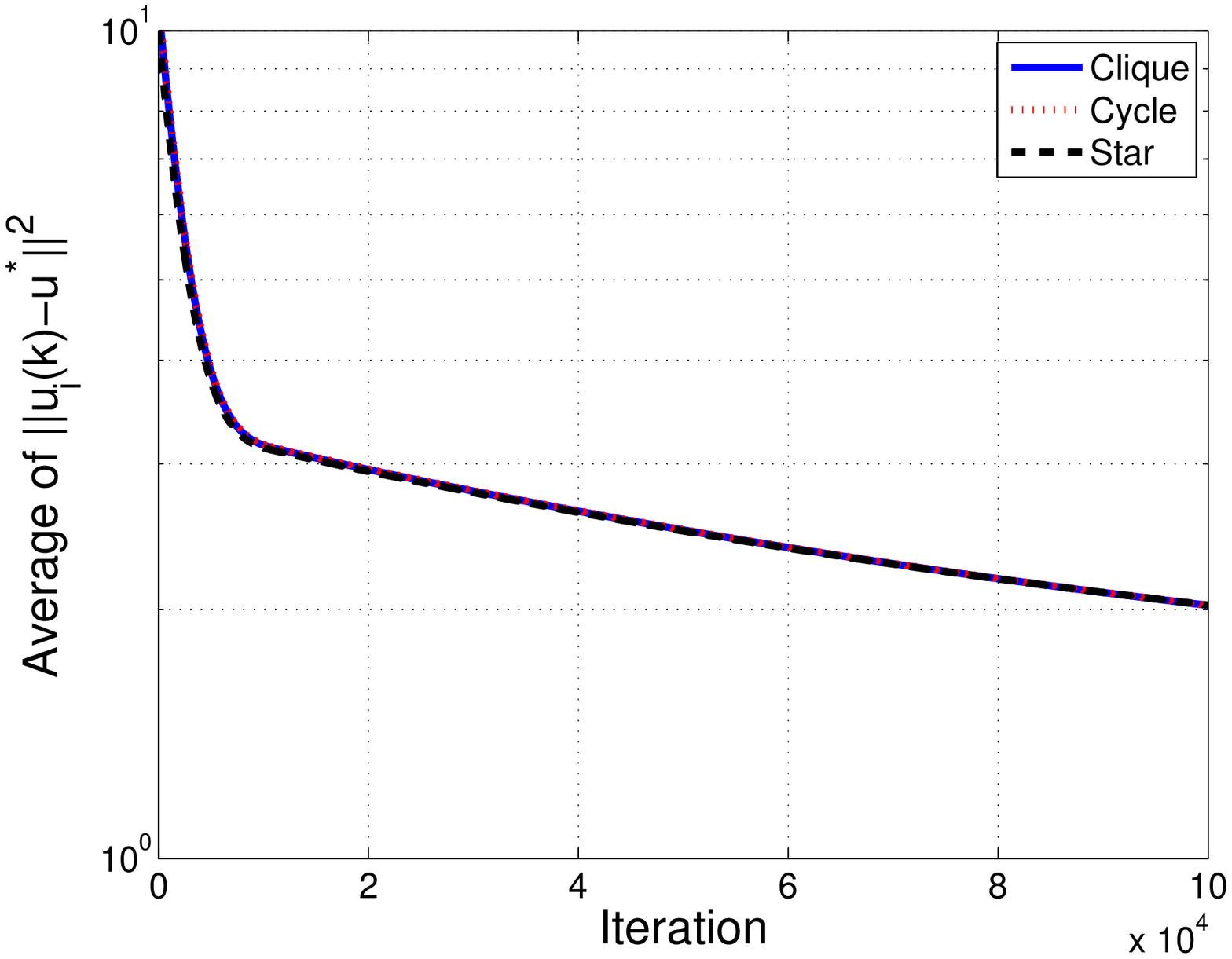}
\end{center}
\caption{\label{fig:m4} Iteration vs $\frac{1}{m}\sum_{i=1}^m\|\mathbf{u}_i(k)-\mathbf{u}^*\|^2$ with a constant stepsize when $m=4$ (left) and $m=10$ (right)}
\end{figure}

In Figures~\ref{fig:m4-dim} and~\ref{fig:m4},
we depict $\frac{1}{m}\sum_{i=1}^m\|\mathbf{u}_i(k)-\mathbf{u}^*\|^2$ over 40,000 and 100,000 iterations when the diminishing and  constant stepsize are used, respectively.
The optimal solution $\mathbf{u}^*$ was obtained by solving the equivalent QP problem (i.e., problem (\ref{eqn:mpc}) with constraints (\ref{eqn:cons1})-(\ref{eqn:cons2}) and (\ref{eqn:qpa})-(\ref{eqn:qpb})) using a commercial QP solver.

We can observe for both cases that the errors go down fast.
An interesting observation is that the network topology does not affect the algorithm performance when the diminishing stepsize is used.
When the constant stepsize is used for the $m=4$ case, star network converges much slower than the other two networks.
This is because the agent selection probability $\g_i$ is different for the center node and the peripheral nodes. As the bound in Proposition \ref{prop:sc} captures, a more aggressive stepsize $\a_i$ should have been used for the peripheral nodes.
For the $m=10$ case, however, the difference is not as clearly visible as in the $m=4$ case. This can be explained by the almost the same spectral gap $1-\sqrt{\lambda}$ (as shown in Proposition~\ref{prop:sc} and Table~\ref{tbl:lambda}).

\begin{table}[t]
\caption{\label{tbl:lambda} Number of agents and $\lambda$}
\begin{center}
\begin{tabular}{|c|c|c|c|}\hline
$m$ & Clique & Cycle & star\\\hline
4 & 0.6667 & 0.7500 & 0.8333\\\hline
10 & 0.8889 & 0.9809 & 0.9444\\\hline
\end{tabular}
\end{center}
\end{table}

\section{Conclusions  \label{sec:con}}
We have considered a distributed problem of minimizing the sum of agents' objective functions over a distributed constraint set $\Xc_i$.
We proposed an asynchronous gossip-based random projection algorithm for solving the problem over a network.
We studied the convergence properties of the algorithm for a random diminishing stepsize and a constant deterministic stepsize.
We established convergence with probability 1~to an optimal solution when the
diminishing stepsizes are used and an error bound when constant stepsizes are used.
We have also provided a simulation result for a distributed robust model predictive control problem.

\bibliographystyle{ieeetran}
\bibliography{soomin002}

\appendix
\subsection{Proof of Lemma \ref{lem:second}\label{app:lemsecond}}
We begin with a lemma which provides some basic relations for
a vector $\check{x}\in \mathcal{Y}$, an arbitrary point $z\in\mathbb{R}^d$,
and two consecutive iterates $x$ and $y$ of a projected-gradient algorithm. The auxiliary point $z$ will be
used to accommodate the iterations $v_i(k)$ of the GRP method which may not belong to the constraint set
$\mathcal{X}$, while $\check{x}$ will be a suitably chosen point in $\mathcal{X}$.

\begin{lemma}\label{lem:first}
Let $\mathcal{Y} \subseteq \mathbb{R}^d$ be a closed convex set.
Let the function $\phi:\mathbb{R}^d\to\mathbb{R}$ be convex and
differentiable over  $\mathbb{R}^d$ with Lipschitz continuous gradients with a constant $L$.
Let $y$ be given by $y = \mathsf{\Pi}_{\mathcal{Y}}[x-\alpha\nabla \phi(x)]$
for some $x \in \mathbb{R}^d \text{ and } \alpha > 0.$
Then,  for all $\check{x} \in \mathcal{Y}$ and $z \in \mathbb{R}^d$, we have:
\begin{itemize}
\item[(a)] For any scalars $\tau_1,\tau_2 >0$,
\begin{align*}
\|y-\check{x}\|^2 \le {}
& (1+8\a^2L^2)\|x-\check{x}\|^2 - 2\a \left(\phi(z)-\phi(\check{x})\right)
-\frac{3}{4}\|y-x\|^2 \cr
&+ (8+\tau_2)\a^2\|\nabla \phi(\check{x})\|^2
+ \tau_1\a^2L^2\|z-\check{x}\|^2+\left(\frac{1}{\tau_1} +\frac{1}{\tau_2}\right)\|x-z\|^2.
\end{align*}
\item[(b)] In addition, if $\phi$ is strongly convex on $\mathbb{R}^d$
with a constant $\sigma>0$, then for any 
$\tau_1,\tau_2>0$,
\begin{align*}
\|y-\check{x}\|^2 \le {}
& (1-\a\sigma+8\a^2L^2)\|x-\check{x}\|^2 - 2\a \langle\nabla\phi(\check{x}), z -\check{x}\rangle
-\frac{3}{4}\|y-x\|^2 \cr
&+ (8+\tau_2)\a^2\|\nabla \phi(\check{x})\|^2
+\tau_1\a^2L^2\|z-\check{x}\|^2 +\left(\frac{1}{\tau_1} +\frac{1}{\tau_2}\right)\|x-z\|^2.
\end{align*}
\end{itemize}
\end{lemma}

\begin{IEEEproof}
For part (a), from the relation defining $y$ and
the strictly non-expansive projection property in~\eqref{eq:proj}, we obtain for any $\check{x} \in \Yc$,
\begin{align}\label{eq:1}
\|y-\check{x}\|^2 \le {}& \|x-\check{x}\|^2 - 2\a\langle\nabla \phi(x),x-\check{x}\rangle
-\|y-x\|^2 + 2\a\langle\nabla \phi(x), x-y\rangle.
\end{align}
We next estimate the term $2\a\langle\nabla \phi(x), x-y\rangle.$
By using Cauchy-Swartz inequality we obtain
$2\a\langle\nabla \phi(x), x-y\rangle \le 2 \alpha\|\nabla \phi(x)\| \|x-y\|.$
By writing $2 \alpha\|\nabla \phi(x)\| \|x-y\|=2 (2\alpha\|\nabla \phi(x)\|) (\|x-y\|/2)$, we find that
\begin{align}\label{eq:2}
2\a\langle\nabla \phi(x), x-y\rangle \le 4\a^2\|\nabla \phi(x)\|^2+\frac{1}{4}\|x-y\|^2.\end{align}
Furthermore, we have $\|\nabla \phi(x)\|^2\le \|(\nabla \phi(x)-\nabla \phi(\check{x}))
+\nabla \phi(\check{x})\|^2$, which by the square-function property $(a+b)^2\le 2(a^2 +b^2)$ yields
$\|\nabla \phi(x)\|^2\le 2\|\nabla \phi(x)-\nabla \phi(\check{x})\|^2 + 2\|\nabla \phi(\check{x})\|^2.$
The preceding relation and the Lipschitz gradient property of $\phi$ imply
\begin{align}\label{eq:3}
\|\nabla \phi(x)\|^2\le 2L\|x-\check{x}\|^2 + 2\|\nabla \phi(\check{x})\|^2.\end{align}
Therefore, from~\eqref{eq:1}--\eqref{eq:3} we obtain
\begin{align}\label{eq:4}
\|y-\check{x}\|^2 \le {}& (1+8\a^2L^2)\|x-\check{x}\|^2 - 2\a\langle\nabla \phi(x),x-\check{x}\rangle
-\frac{3}{4}\|y-x\|^2 + 8\a^2\|\nabla \phi(\check{x})\|^2.
\end{align}
Next, we estimate the term $2\a\langle\nabla \phi(x),x-\check{x}\rangle$  using the convexity of $\phi$,
\begin{align}\label{eq:5}
\langle\nabla \phi(x),x-\check{x}\rangle
\ge \phi(x)-\phi(\check{x}) = \left(\phi(x)-\phi(z) \right)+\left(\phi(z)-\phi(\check{x})\right),
\end{align}
where $z\in\mathbb{R}^d$ is some given point. It remains to bound the term $\phi(x)-\phi(z)$, for which by
convexity of $\phi$ we further have
\begin{align*}
\phi(x)-\phi(z)\ge \langle\nabla\phi(z),x-z\rangle\ge -\|\nabla \phi(z)\|\,\|x-z\|.
\end{align*}
By writing $\|\nabla \phi(z)\|\le \|\nabla \phi(z) - \nabla\phi(\check{x})\| +\|\nabla\phi(\check{x})\|$
and using the Lipschitz-gradient property of $\phi$, we obtain
\begin{align*}
\phi(x)-\phi(z)\ge -L\|z-\check{x}\|\,\|x-z\|-\|\nabla\phi(\check{x})\|\,\|x-z\|.
\end{align*}
Multiplying the preceding relation with $2\a$ and using
$2\a L\|z-\check{x}\|\,\|x-z\|=2 (\a\sqrt{\tau_1} L\|z-\check{x}\|)(\|x-z\|/\sqrt{\tau_1})\le
\tau_1\a^2 L^2\|z-\check{x}\|^2 +\|x-z\|^2/\tau_1$,
$2\a \|\nabla\phi(\check{x})\|\,\|x-z\|=2(\a \sqrt{\tau_2}\|\nabla\phi(\check{x})\|)
(\|x-z\|/\sqrt{\tau_2})\le
\tau_2\a^2\|\nabla\phi(\check{x})\|^2 +\|x-z\|^2/\tau_2$ for some $\tau_1,\tau_2>0$, we obtain
\begin{align}\label{eq:6}
2\a \left(\phi(x)-\phi(z)\right)\ge
-\tau_1\a^2 L^2\|z-\check{x}\|^2 - \tau_2\a^2\|\nabla\phi(\check{x})\|^2 -
\left(\frac{1}{\tau_1} +\frac{1}{\tau_2}\right) \|x-z\|^2.
\end{align}
Thus, from Eqs.~\eqref{eq:4}--\eqref{eq:6} it follows that
\begin{align}\label{eq:7}
\|y-\check{x}\|^2 \le {}
& (1+8\a^2L^2)\|x-\check{x}\|^2 - 2\a \left(\phi(z)-\phi(\check{x})\right)
-\frac{3}{4}\|y-x\|^2 \cr
&+ (8+\tau_2)\a^2\|\nabla \phi(\check{x})\|^2
+\tau_1\a^2L^2\|z-\check{x}\|^2 + \left(\frac{1}{\tau_1} +\frac{1}{\tau_2}\right)\|x-z\|^2,
\end{align}
thus proving the relation in part (a).
The relation in part (b) follows similarly by using the strong convexity of $\phi$ in Eq.~\eqref{eq:5},
i.e., $\langle\nabla \phi(x),x-\check{x}\rangle\ge
\phi(x)-\phi(\check{x})+\frac{\sigma}{2}\|x-\check{x}\|^2$ for all $x,\check{x}\in\mathbb{R}^d$.
\end{IEEEproof}

The proof of Lemma~\ref{lem:second} relies on Lemma \ref{lem:first}(a) and
the fact that the event $E_i(k) = \{i \in \{I_k,J_k\}\}$ that agent $i$ updates at any time is
independent of the past. Due to this, the number of updates that any agent $i$ has performed
until time $k$ behaves almost as $1/k$ when $k$ is large enough.
The long term estimates for the stepsize $\a_i(k) = \frac{1}{\Gamma_i(k)}$ in terms of the probability
$\g_i$ that agent $i$ updates are given in the following lemma.

\begin{lemma}\label{lem:stepsize}(see \cite{Nedic11})
Let $\a_i(k) = 1/\Gamma_i(k)$ for all $k \ge 1$ and $i \in V$.
Let $\pi_{\min} = \min_{\{i,j\}\in E} \pi_{ij}$.
Also, let $q$ be a constant such that $0 < q <1/2$.
Then, there exists a large enough $\tilde{k}$ (which depends on $q$ and $m$) such that
with probability 1 for all $k \ge \tilde{k}$ and $i \in V$,
\begin{enumerate}
\item[(a)] $\displaystyle \a_i(k) \le \frac{2}{k\g_i}$,\quad
(b) $\displaystyle \a_i^2(k) \le \frac{4m^2}{k^2(1+\pi_{\min})^2}$,
\item[(c)] $\displaystyle \left|\a_i(k)-\frac{1}{k\g_i}\right| \le \frac{2}{k^{\frac{3}{2}-q}(1+\pi_{\min})^2}$.
\end{enumerate}
\end{lemma}
According to this lemma, the stepsizes $\a_i(k)$ exhibit the same behavior
as the deterministic stepsize $1/k$ in a long run. The result is critical for dealing with
the cross dependencies of the random stepsizes and the other randomness in the GRP method.

{\bf Proof of Lemma~\ref{lem:second}.} Consider $i \in \{I_k, J_k\}$, and
use Lemma \ref{lem:first}(a) with the following identification:
$\Yc = \Xc_i^{\Omega_i(k)}$ and $\check{x}\in\Xc\subseteq\Xc_i^{\Omega_i(k)}$,
$y = x_i(k)$, $x = v_i(k)$, $z = z_i(k) \triangleq \proj_{\Xc}[v_i(k)]$, $\phi = f_i$, and
$\a = \a_i(k)$.
Then, for any $\check{x} \in \Xc$ and $k \ge 1$,
\begin{align*}
&\|x_i(k)-\check{x}\|^2 \le {}
 (1+8\a^2_i(k) L_i^2)\|v_i(k)-\check{x}\|^2 - 2\a_i(k) \left(f_i(z_i(k))-f_i(\check{x})\right)
-\frac{3}{4}\|x_i(k)-v_i(k)\|^2 \cr
&+ (8+\tau_2)\a^2_i(k)\|\nabla f_i(\check{x})\|^2
+\tau_1\a^2_i(k)L_i^2\|z_i(k)-\check{x}\|^2 + \left(\frac{1}{\tau_1} +\frac{1}{\tau_2}\right)\|v_i(k)-z_i(k)\|^2,
\end{align*}
By Assumption~\ref{assume:f}(d), we have $\|\nabla f_i(\check{x})\|\le G_f$.
Further, we let $\tau_1=\tau_2=4\eta$ for some $\eta>0$,
 and by using Lemma~\ref{lem:stepsize}(b) we find that w.p.1 for all $k$ large enough,
\begin{align}\label{eq:8}
\|x_i(k)-\check{x}\|^2 \le {} &
 \left(1+\frac{c_1}{k^2}\right)\|v_i(k)-\check{x}\|^2 - 2\a_i(k) \left(f_i(z_i(k))-f_i(\check{x})\right)\cr
& -\frac{3}{4}\|x_i(k)-v_i(k)\|^2
+ \frac{c_2}{k^2}
+\frac{c_3}{k^2}\|z_i(k)-\check{x}\|^2 +\frac{1}{2\eta}\|v_i(k)-z_i(k)\|^2,
\end{align}
where $c_1 = \frac{32m^2 \bar L^2}{(1+\pi_{\min})^2}$, $\bar L=\max_{i}L_i$,  $c_2=\frac{4(8+4\eta)m^2G_f^2}{(1+\pi_{\min})^2}$
and $c_3=\frac{16\eta m^2\bar L^2}{(1+\pi_{\min})^2}$.
Next, consider  $2\a_i(k) \left(f_i(z_i(k))-f_i(\check{x})\right)$, for which we can write
\[2\a_i(k) \left(f_i(z_i(k))-f_i(\check{x})\right)
\ge\frac{2}{k\g_i} \left(f_i(z_i(k))-f_i(\check{x})\right)  -2\left|\a_i(k)-\frac{1}{k\g_i}\right|\,
\left|f_i(z_i(k))-f_i(\check{x})\right|.\]
Since $f_i$ has bounded gradients over the set $\Xc$, it is Lipschitz continuous over $\Xc$.
Thus, since $z_i(k),\check{x}\in \Xc$, it follows that
$\left|f_i(z_i(k))-f_i(\check{x})\right|\ge G_f\|z_i(k)-\check{x}\|$. This and Lemma~\ref{lem:stepsize}$($c)
imply
\begin{align*}
2\a_i(k) \left(f_i(z_i(k))-f_i(\check{x})\right)
&\ge \frac{2}{k\g_i} \left(f_i(z_i(k))-f_i(\check{x})\right)
-2\frac{2}{k^{\frac{3}{2}-q}(1+\pi_{\min})^2}G_f\|z_i(k)-\check{x}\|\cr
&\ge \frac{2}{k\g_i} \left(f_i(z_i(k))-f_i(\check{x})\right)
-\frac{2}{k^{\frac{3}{2}-q}(1+\pi_{\min})^2}\left(G_f^2 + \|z_i(k)-\check{x}\|^2\right),
\end{align*}
where the last inequality follows by the Cauchy-Schwarz inequality.
Combining the preceding relation with Eq.~\eqref{eq:8}, we obtain w.p.1 for $k$ large enough
\begin{eqnarray}\label{eq:9}
\|x_i(k)-\check{x}\|^2 &\le  &
 \left(1+\frac{c_1}{k^2}\right)\|v_i(k)-\check{x}\|^2 - \frac{2}{\g_ik} \left(f_i(z_i(k))-f_i(\check{x})\right)
 -\frac{3}{4}\|x_i(k)-v_i(k)\|^2 \cr
&&+ \frac{c_2}{k^2} +\frac{c_4}{k^{\frac{3}{2}-q} }
+\left(\frac{c_3}{k^2}+\frac{c_5}{k^{\frac{3}{2}-q}}\right)\|z_i(k)-\check{x}\|^2
+\frac{1}{2\eta}\|v_i(k)-z_i(k)\|^2,
\end{eqnarray}
where $c_4=\frac{2}{(1+\pi_{\min})^2}G_f^2$ and $c_5=\frac{2}{(1+\pi_{\min})^2}$.

By the definition of the projection, we have $\|v_i(k)-z_i(k)\| = \dist(v_i(k),\Xc),$ and
\[
\|x_i(k)-v_i(k)\|\ge \left\|\mathsf{\Pi}_{\mathcal{X}_i^{\Omega_i(k) }}[v_i(k)]-v_i(k)\right\|
=\dist(v_i(k),\mathcal{X}_i^{\Omega_i(k)}).\]
Taking the expectation in (\ref{eq:9}) conditioned jointly on $\Fc_{k-1}, I_k$ and $J_k$,
we obtain for any $\check{x} \in \Xc$, $i \in \{I_k,J_k\}$ w.p.1 for all $k$ large enough
we can see that
\begin{eqnarray*}
&&\EXP{\|x_i(k)-\check{x}\|^2\mid\Fc_{k-1}, I_k,J_k} \le
 \left(1+\frac{c_1}{k^2}\right)\|v_i(k)-\check{x}\|^2 - \frac{2}{\g_ik} \left(f_i(z_i(k))-f_i(\check{x})\right)\cr
&& -\frac{3}{4}\EXP{ \dist^2(v_i(k),\mathcal{X}_i^{\Omega_i(k)})\mid v_i(k)}
+ \frac{c_6}{k^{\frac{3}{2}-q} }
+\frac{c_7}{k^{\frac{3}{2}-q}}\|z_i(k)-\check{x}\|^2
+\frac{1}{2\eta}\dist^2(v_i(k),\Xc),
\end{eqnarray*}
where $c_6=c_2+c_4$ and $c_7=c_3+c_5$.
Using the regularity condition (Assumption~\ref{assume:c}), we have
\[
\EXP{\dist^2(v_i(k),\mathcal{X}_i^{\Omega_i(k)})\mid v_i(k)}\ge
\frac{1}{c}\dist^2(v_i(k),\Xc).
\]
Thus, by letting $\eta=c$,
from the preceding two relations we have w.p.1 for all $k$ large enough
\begin{eqnarray*}
&&\EXP{\|x_i(k)-\check{x}\|^2\mid\Fc_{k-1}, I_k,J_k} \le
 \left(1+\frac{a_1}{k^2}\right)\|v_i(k)-\check{x}\|^2 - \frac{2}{\g_ik} \left(f_i(z_i(k))-f_i(\check{x})\right)\cr
&& -\frac{1}{4c} \dist^2(v_i(k),\Xc)
+ \frac{c_6}{k^{\frac{3}{2}-q} }
+\frac{c_7}{k^{\frac{3}{2}-q}}\|z_i(k)-\check{x}\|^2.
\end{eqnarray*}
The preceding inequality holds with probability $\g_i$
(when agent $i$ updates), and otherwise
$x_i(k) = v_i(k)$ with probability $1-\g_i$ (when agent $i$ does not update).
Hence, w.p.1 for any $\check{x} \in \Xc$, all $i \in V$, and all $k$ large enough we have
\begin{eqnarray*}
&&\EXP{\|x_i(k)-\check{x}\|^2\mid\Fc_{k-1}} \le
 \left(1+\frac{\g_ia_1}{k^2}\right)\EXP{\|v_i(k)-\check{x}\|^2 \mid\Fc_{k-1}}
 - \frac{2}{k} \EXP{f_i(z_i(k))-f_i(\check{x})\mid\Fc_{k-1}}\cr
&& -\frac{\g_i}{4c} \EXP{\dist^2(v_i(k),\Xc)\mid\Fc_{k-1}}
+ \frac{\g_i c_6}{k^{\frac{3}{2}-q} }
+ \frac{\g_ic_7}{k^{\frac{3}{2}-q}} \EXP{\|z_i(k)-\check{x}\|^2\mid\Fc_{k-1}}.
\end{eqnarray*}
Since $\g_i\le 1$, the relation of Lemma~\ref{lem:second} follows
by letting $a_1=c_1,$ $ a_2=c_6$ and $a_3=c_7$.
$\blacksquare$
\subsection{Proof of Lemma \ref{lemma:key}\label{app:lemkey}}
The proof of this lemma and the proofs of the other lemmas, often,
rely on the relations implied by the convexity of the squared-norm. In particular,
by the definition of $v_i(k) $ in~(\ref{eqn:algoa}), the convexity of the squared-norm function and
the doubly stochastic weights $W(k)$, we have for any $x\in\mathbb{R}^d$,
\begin{equation}\label{eqn:ds}
\sum_{i=1}^m \Es[\|v_i(k)-x\|^2 \mid \Fc_{k-1}]
\le \sum_{i=1}^m\sum_{j=1}^m \bar{W}_{ij}\|x_j(k-1)-x\|^2 = \sum_{j=1}^m \|x_j(k-1)-x\|^2.
\end{equation}
Similarly, by the convexity of the distance function $x\mapsto\dist^2(x,\Xc)$ (see \cite[p.~88]{BNO}),
we have
\begin{align}\label{eqn:ds-dist}
\sum_{i=1}^m\Es\left[\dist^2(v_i(k),\Xc) \mid \Fc_{k-1}\right] \le
\sum_{j=1}^m \dist^2(x_j(k-1),\Xc).
\end{align}

{\bf Proof of Lemma \ref{lemma:key}.}
To prove part (a), we start with
Lemma~\ref{lem:second}, where we let $\check{x}= z_i(k)=\Ps_{\Xc}[v_i(k)]$.
Then, for all $k$ large enough and all $i \in V$, we obtain w.p.1,
\begin{eqnarray}\label{eq:oneofmany}
\EXP{\|x_i(k)-\Ps_{\Xc}[v_i(k)]\|^2\mid\Fc_{k-1}} &\le&
 \left(1+\frac{a_1}{k^2}\right)\EXP{\dist^2(v_i(k),\Xc) \mid\Fc_{k-1}}\cr
&&  -\frac{\g_i}{4c} \EXP{\dist^2(v_i(k),\Xc)\mid\Fc_{k-1}}
+ \frac{a_2}{k^{\frac{3}{2}-q} },
\end{eqnarray}
where $q\in(0,1/2)$.
By the definition of the projection, we have
$\dist(x_i(k),\Xc) \le \|x_i(k)-\mathsf{\Pi}_{\mathcal{X}}[v_i(k)]\|.$
Using this relation in Eq.~\eqref{eq:oneofmany} and, then, summing the resulting  relations over $i$ and applying Eq.~\eqref{eqn:ds-dist}, we find that w.p.1 for all $k$ large enough and all $i \in V$,
\begin{eqnarray}\label{eq:extra}
&&\sum_{i=1}^m\EXP{\dist^2(x_i(k),\Xc)\mid\Fc_{k-1}}
\le  \left(1+\frac{a_1}{k^2}\right) \sum_{j=1}^m\dist^2(x_j(k-1),\Xc) \cr
&& -\frac{\underline{\g} }{4c}\sum_{i=1}^m\EXP{\dist^2(v_i(k),\Xc)\mid\Fc_{k-1}}
+ \frac{a_2m}{k^{\frac{3}{2}-q} },
\end{eqnarray}
where $\underline{\g} = \min_i \g_i$.
Therefore, for all $k$ large enough, the conditions of Lemma \ref{thm:super} are satisfied
(for a time-delayed sequence), so we conclude that
$\sum_{k=1}^\infty \EXP{\dist^2(v_i(k),\Xc)\mid\Fc_{k-1}}<\infty$ for all~$i$.
Taking the total expectation in relation (\ref{eq:extra}), it also follows that $\sum_{k=1}^\infty \EXP{\dist^2(v_i(k),\Xc)}<\infty$
for all $i \in V$, which by the Monotone Convergence Theorem \cite[p.92]{Royden-book}
yields $\lim_{k\to\infty} \dist(v_i(k), \Xc)=0$ for all $i$ w.p.1, showing the result in part (a).


For part (b), note that for $\|e_i(k)\|$, using $z_i(k)=\mathsf{\Pi}_{\mathcal{X}}[v_i(k)]$,
we can write for $i \in \{I_k,J_k\}$,
\begin{align*}
\|e_i(k)\| &\le\|x_i(k) - z_i(k)\| +\|z_i(k) - v_i(k)\|\cr
&=\left\|\mathsf{\Pi}_{\mathcal{X}_i^{\Omega_i(k)}}[v_i(k) - \alpha_i(k) \nabla f_i(v_i(k))]
- z_i(k)\right\|+\left\|z_i(k)- v_i(k)\right\|.
\end{align*}
Since $\mathcal{X} \subseteq \mathcal{X}_i^{\Omega_i(k)}$ and
$z_i(k) \in \mathcal{X}$, we have $z_i(k)\in \mathcal{X}_i^{\Omega_i(k)}$.
Using the projection non expansiveness property of Eq.~\eqref{eq:proj-c}, we obtain
\begin{align*}
\|e_i(k)\|
&\leq  \|v_i(k) - \alpha_i(k) \nabla f_i(v_i(k)) -z_i(k)\|+\|z_i(k)- v_i(k)\|\nonumber\\
&\leq  2\| v_i(k)- z_i(k)\| + \alpha_i(k) \|\nabla f_i(v_i(k))\|.
\end{align*}
Further, from the Lipschitz gradient property of $f_i$
and the gradient boundedness property (Assumptions~\ref{assume:f}(c) and \ref{assume:f}(d))
it follows that
\begin{align}\label{eqn:e_final}
\|e_i(k)\|
&\le 2\|v_i(k) - z_i(k)\|
+ \alpha_i(k) \left(\|\nabla f_i(v_i(k))-\nabla f_i(z_i(k))\| +\|\nabla f_i(z_i(k))\|\right)\nonumber\\
&\le \left(2+\a_i(1) L_i\right) \dist (v_i(k),\Xc)+ \alpha_i(k) G_f,
\end{align}
where the last inequality follows by $\a_i(k)\le \a_i(1)$ and $\|v_i(k) - z_i(k)\| = \dist (v_i(k),\Xc)$.
Using the Cauchy-Schwartz inequality and Lemma \ref{lem:stepsize}(a) (i.e.
$\a_i(k)\le 2/(k\g_i)$), we have for all $i\in \{I_k,J_k\}$ and $k\ge \tilde{k}$,
\begin{equation}\label{eqn:e}
\|e_i(k)\|^2 \leq 2(2+\a_i(1) L_i)^2\dist^2 (v_i(k),\Xc) + \frac{8m^2}{k^2} G_f^2,
\end{equation}
where we also use $\g_i\ge\frac{1}{m}$.
Taking the expectation in (\ref{eqn:e}) conditioned on $\Fc_{k-1}, I_k, J_k$
and noting that the preceding inequality holds with probability
$\g_i$, and $x_i(k) = v_i(k)$ with probability $1-\g_i$, we obtain with probability 1 for all
$k \ge \tilde{k}$ and $i \in V$,
\begin{equation*}
\Es[\|e_i(k)\|^2\mid \Fc_{k-1}]
\leq 2\g_i(2+\a_i(1) \bar L)^2\Es[\dist^2(v_i(k),\Xc)\mid \Fc_{k-1}] + \frac{8\g_im^2}{k^2} G_f^2,
\end{equation*}
where $\bar L=\max_i L_i$.
By part (a) of this lemma, we have $\sum_{k=1}^\infty \Es[\dist^2(v_i(k),\Xc)\mid \Fc_{k-1}] <\infty$
w.p.1 for all $i$.
As $\sum_{k=1}^\infty\frac{1}{k^2}<\infty$, we conclude that
$\sum_{k=1}^\infty \Es[\|e_i(k)\|^2\mid \Fc_{k-1}] <\infty$ for all $i \in V$ w.p.1.
Furthermore, by relation~\eqref{eqn:e} and part (a) of the lemma we find that
$\lim_{k\to\infty} \|e_i(k)\|=0$ for all $i$ w.p.1.
\quad $\blacksquare$
\subsection{Proof of Lemma \ref{lem:disagree}\label{app:lemdisagree}}
The proof of this Lemma makes use of an additional result, which is given below.

\begin{lemma}\label{lem:coordconsensus}
Let $\{W(k)\}$ be an iid sequence of $m\times m$ symmetric and stochastic matrices.
Consider a sequence $\{\t(k)\}\subset \mathbb{R}^m$ generated by the following dynamics
\begin{equation}\label{eqn:conerr}
\t(k) = W(k) \t(k-1) + \e(k)\qquad\hbox{for }k\ge1.
\end{equation}
Then, we have with probability 1 for all $k\ge1$,
\[
\Es[\|\Delta(k)\| \mid \Fc_{k-1}] \le \sqrt{\lambda}\|\Delta(k-1)\| + \Es[\|\e(k)\| \mid \Fc_{k-1}].
\]
where $\Delta(k) \triangleq \t(k)-\frac{1}{m}\1b\1b^T\t(k)$
and $\lambda < 1$ is the second largest eigenvalue of $\bar{W}=\EXP{W(k)}$.
\end{lemma}
\begin{IEEEproof}
Define the sequences of averaged coordinate values as $\t_k^{\ave} \triangleq \frac{1}{m}\sum_{i=1}^m \t_i(k)$
and $\e_k^{\ave} \triangleq \frac{1}{m}\sum_{i=1}^m \e_i(k)$.
From relation \eqref{eqn:conerr}, by taking averages over the coordinates, we have
$\t_k^{\ave} = \t_{k-1}^{\ave} + \e_k^{\ave}.$
Using $\1b \in \mathbb{R}^m$, a vector with all its elements 1, we can write
$\t_k^{\ave}\1b = \t_{k-1}^{\ave}\1b + \e_k^{\ave}\1b,$
or equivalently,
\begin{equation}\label{eqn:avgvec}
\frac{1}{m}\1b\1b^T\t(k) = \frac{1}{m}\1b\1b^T \t(k-1) + \frac{1}{m}\1b\1b^T\e(k).
\end{equation}
From equations \eqref{eqn:conerr} and \eqref{eqn:avgvec}, we have:
\[
\t(k)-\frac{1}{m}\1b\1b^T\t(k) = \left(W(k)-\frac{1}{m}\1b\1b^T\right) \t(k-1) + \left(I-\frac{1}{m}\1b\1b^T\right)\e(k).
\]
Since $\left(W(k)-\frac{1}{m}\1b\1b^T\right)\frac{1}{m}\1b\1b^T\t(k-1)
= \left(W(k)-\frac{1}{m}\1b\1b^T\right)\t_{k-1}^{\ave}\1b = 0$,
by letting $\Delta(k) \triangleq \t(k)-\frac{1}{m}\1b\1b^T\t(k)$,
$D_k \triangleq W(k)-\frac{1}{m}\1b\1b^T$ and $M \triangleq I-\frac{1}{m}\1b\1b^T$,
 it follows that
$\Delta(k) = D_k \Delta(k-1) + M\e(k)$ for all $k\ge1$.
By taking the norm and the expectation conditioned on the history $\Fc_{k-1}$,
from the preceding relation we have w.p.1 for $k \ge 1$,
\begin{equation}\label{eqn:rec}
\Es[\|\Delta(k)\| \mid \Fc_{k-1}] = \Es[\|D_k \Delta(k-1)\|\mid \Fc_{k-1}] + \Es[\|M\e(k)\| \mid \Fc_{k-1}].
\end{equation}

From Eq.~\eqref{eq:ly} and the fact that $W(k)$ is independent of the past $\Fc_{k-1}$, we obtain
$\Es[\|D_k\Delta(k-1)\|^2\mid \Fc_{k-1}] 
\le \lambda \|\Delta(k-1)\|^2,$
where $\lambda$ is the second largest eigenvalue of the matrix $\bar{W}$.
Using $\Es[\|x\|] \le \sqrt{\Es[\|x\|^2]}$, we obtain
$\Es[\|D_k\Delta(k-1)\|\mid \Fc_{k-1}] \le \sqrt{\lambda} \|\Delta(k-1)\|$ for all $k \ge 1$.
For the second term on the right hand side of \eqref{eqn:rec}, we have
$\Es[\|M\e(k)\| \mid \Fc_{k-1}] = \Es[\|\e(k)\| \mid \Fc_{k-1}],$
since $M = I-\frac{1}{m}\1b\1b^T$ is a projection matrix and, thus, $\|M\| = 1$.
\end{IEEEproof}

{\it Proof of Lemma \ref{lem:disagree}.}
We consider coordinate-wise relations by defining the vector $y_{\ell}(k) \in \mathbb{R}^m$ for
$\ell = 1,\ldots,d$ such that $[y_{\ell}(k)]_i = [x_i(k)]_{\ell}$ for all $i$.
From algorithm (\ref{eqn:algoa})-(\ref{eqn:algoc}), we have
\begin{equation*}
y_{\ell}(k) = W(k)y_{\ell}(k-1) + \delta_{\ell}(k) \quad \textrm{for } k \ge 1,
\end{equation*}
where $\delta_{\ell}(k) \in \mathbb{R}^m$ is a vector whose coordinates are defined by
\begin{equation}\label{eqn:delta}
[\delta_{\ell}(k)]_i =
\left[\mathsf{\Pi}_{\mathcal{X}_i^{\Omega_i(k)}}[v_{i}(k)-\alpha_i(k)\nabla f(v_{i}(k))]-v_{i}(k)\right]_{\ell}
\qquad\hbox{if  $i \in \{I_k, J_k\},$}
\end{equation}
and otherwise $[\delta_{\ell}(k)]_i =0$. 
Since the matrices $W(k)$ are doubly stochastic for all $k\ge 1$,
from Lemma \ref{lem:coordconsensus} we obtain
\begin{equation}\label{eqn:eqrec}
\Es[\|y_{\ell}(k)-[\bar{x}(k)]_{\ell}\1b\| \mid \Fc_{k-1}] \le \sqrt{\lambda}\|y_{\ell}(k-1)-[\bar{x}(k-1)]_{\ell}\1b\| + \Es[\|\delta_{\ell}(k)\| \mid \Fc_{k-1}],
\end{equation}
where $[\bar{x}(k)]_{\ell} = \frac{1}{m}\1b^Ty_{\ell}(k)$ and $\lambda < 1$ by Assumption~\ref{assume:1}.

We next consider $\delta_{\ell}(k)$ as given by~(\ref{eqn:delta}), for which we have for all $k \ge 1$,
\[\|\delta_{\ell}(k)\|^2 \le
\sum_{i=1}^m\left\|\mathsf{\Pi}_{\mathcal{X}_i^{\Omega_i(k)}}[v_{i}(k)-\alpha_i(k)(\nabla f_i(v_{i}(k)))]-v_{i}(k)\right\|^2.\]
Letting $z_i(k) \triangleq \mathsf{\Pi}_{\Xc}[v_i(k)]$, observing that $z_i(k) \in \Xc_i^{\Omega_i(k)}$
and using the projection property in Eq.~\eqref{eq:proj-c},  we obtain
\begin{align*}
\|\delta_{\ell}(k)\|^2
{}& \le
\sum_{i=1}^m\left(\left\|\mathsf{\Pi}_{\mathcal{X}_i^{\Omega_i(k)}}[v_{i}(k)-\alpha_i(k)\nabla f_i(v_{i}(k))]-z_i(k)\right\| + \|z_i(k)-v_{i}(k)\|\right)^2\cr
{}&\le \sum_{i=1}^m\left(\alpha_i(k)\left\|\nabla f_i(v_{i}(k))\right\| + 2\|z_i(k)-v_{i}(k)\|\right)^2 .
\end{align*}
Applying the Cauchy-Schwartz inequality, we can obtain
\[\|\delta_{\ell}(k)\|^2
\le \sum_{i=1}^m\left(2\alpha_i^2(k)\left\|\nabla f_i(v_{i}(k))\right\|^2 + 4\|z_i(k)-v_{i}(k)\|^2\right).\]
The term $\|\nabla f_i(v_i(k))\|^2$ can be further evaluated by using the Lipschitz property
and the bounded gradient assumption (Assumption \ref{assume:f}(d)),
\[
\|\nabla f_i(v_i(k))\|^2 \le 2\|\nabla f_i(v_i(k))-\nabla f_i(z_i(k))\|^2 + 2\|\nabla f_i(z_i(k))\|^2
\le  2L^2_i\|v_i(k)-z_i(k)\|^2 + 2G_f^2.
\]
From Lemma \ref{lem:stepsize}(b), there exists a large enough $\tilde{k}$ such that
$\a_i^2(k) \le 4m^2/k^2\le 4m^2/\tilde k^2$ w.p.1 for all $k \ge \tilde{k}$.
Therefore, noting that $\|z_i(k)-v_{i}(k)\|=\dist(v_i(k),\Xc)$,
we obtain  for all $k \ge \tilde{k}$ with probability 1,
\begin{align*}
\|\delta_{\ell}(k)\|^2  {}& \le
\left(4+\frac{16m^2}{\tilde{k}^2}\bar{L}^2\right)\sum_{i=1}^m\dist^2(v_i(k),\Xc)+ \frac{16m^2}{k^2}G_f^2,
\end{align*}
with $\bar L=\max_i L_i$.
Taking the expectation with respect to $\Fc_{k-1}$ and using
$\Es[\|x\|] \le \sqrt{\Es[\|x\|^2]}$, we obtain $\Es[\|\delta_{\ell}(k)\| \mid \Fc_{k-1}] \le b_k$,
where
\[b_k=
\sqrt{ \left(4+\frac{16m^2}{\tilde{k}^2}\bar L^2\right)\sum_{i=1}^m\Es[\dist^2(v_i(k),\Xc)\mid \Fc_{k-1}]+ \frac{16m^2}{k^2}G_f^2}.\]
From the preceding and
relation~(\ref{eqn:eqrec}), we obtain for all $k \ge \tilde{k}$ with probability 1,
\begin{align} \label{eqn:eqsuper}
\frac{1}{k}\Es[\|y_{\ell}(k)-[\bar{x}(k)]_{\ell}\1b\| \mid \Fc_{k-1}]
\le {}&  \frac{1}{k-1} \|y_{\ell}(k-1)-[\bar{x}(k-1)]_{\ell}\1b\|\cr
{}&- \frac{1-\sqrt{\lambda}}{k} \|y_{\ell}(k-1)-[\bar{x}(k-1)]_{\ell}\1b\|
 + \frac{1}{k} b_k.
\end{align}
Noting that $\frac{1}{k} b_k\le (1/k^2 + b_k^2)/2$, and that $\sum_{k=1}^\infty b_k^2<\infty$
by Lemma \ref{lemma:key}(a), the term $\frac{1}{k} b_k$ is summable.
From this and the fact that $1-\sqrt{\lambda} >0 $, relation (\ref{eqn:eqsuper}) satisfies all the conditions in Lemma \ref{thm:super}. It follows that
$\sum_{k=1}^{\infty}\frac{1}{k} \Es[\|y_{\ell}(k)-[\bar{x}(k)]_{\ell}\1b\|\mid\Fc_{k-1}]  < \infty$ with probability 1
for any $\ell = 1,\ldots,d$.
This and the definition of $y_{\ell}(k)$  implies that with probability 1
\begin{equation}\label{eqn:mres}
\sum_{k=1}^{\infty}\frac{1}{k} \Es[\|x_i(k)-\bar{x}(k)\|\mid\Fc_{k-1}]  < \infty \textrm{ for all } i \in V,
\end{equation}
where $\bar x(k)=\sum_{j=1}^m x_j(k)$.
Next, consider $\|v_i(k)- \bar v(k)\|$.
Since $v_i(k)=\sum_{j=1}^m [W(k)]_{ij}\, x_j(k-1)$ (see~\eqref{eqn:algoa})
and $W(k)$ is doubly stochastic, by using the
convexity of the norm, for $\bar v(k)=\frac{1}{m}\sum_{j=1}^m v_j(k)$
we can see that
$\sum_{i=1}^m\|v_i(k)- \bar v(k)\|\le \sum_{j=1}^m \left\|x_j(k-1)- \bar x(k-1)\right\|$.
By using relation~\eqref{eqn:mres}, we conclude that
$\sum_{k=1}^\infty\frac{1}{k}\Es[\|v_i(k)- \bar v(k)\|\mid\Fc_{k-1}]<\infty$ for all $i\in V$ w.p.1.
\quad $\blacksquare$
\subsection{Proof of Lemma \ref{lem:base}\label{app:lembase}}
Let  $i\in\{I_k, J_k\}$. Then, using the definition of
the iterate $x_i(k)$ in (\ref{eqn:algoa})-(\ref{eqn:algoc}), and Lemma \ref{lem:first}(b) with the following
identification: $\mathcal{Y} = \mathcal{X}_i^{\Omega_i(k)}$,
$y = x_i(k)$, $x = v_i(k)$, $z = z_i(k) = \mathsf{\Pi}_{\mathcal{X}}[v_i(k)]$,
$\alpha = \alpha_i$, $\check{x} = x\in\Xc$, $\phi = f_i$,
$L=L_i$ and $\tau_1=\tau_2=8c$, we obtain
\begin{align*}
&\|x_i(k)-x\|^2
\leq  \left(1- \sigma_i\alpha_i+8\alpha_i^2L_i^2\right)\|v_i(k)-x\|^2
- 2\alpha_i\langle \nabla f_i(x), z_i(k)-x\rangle\nonumber\\
{}&   - \frac{3}{4}\|x_i(k)-v_i(k)\|^2
+ 8(1+c)\a_i^2\|\nabla f_i(x)\|^2
+4c\a_i^2L_i^2\|z_i(k)-x\|^2+ \frac{1}{4c} \|v_i(k)-z_i(k)\|^2.
\end{align*}
By Assumption~\ref{assume:f}(d), we have $\|\nabla f_i(x)\|\le G_f$, while
by the non-expansiveness projection property we have $\|z_i(k)-x\|\le \|v_i(k)-x\|$.
Furthermore, $\|v_i(k)-z_i(k)\| = \text{dist}(v_i(k),\mathcal{X})$ since
 $z_i(k) = \mathsf{\Pi}_{\mathcal{X}}[v_i(k)]$.
Therefore, for all $k\ge1$ and $i\in\{I_k, J_k\}$,
\begin{align}\label{eq:remove}
&\|x_i(k)-x\|^2
\leq  \left(1- \sigma_i\alpha_i+4(2+c)\alpha_i^2L_i^2\right)\|v_i(k)-x\|^2
- 2\alpha_i\langle \nabla f_i(x), z_i(k)-x\rangle\nonumber\\
{}&  + 8(1+c)\a_i^2G_f^2- \frac{3}{4}\|x_i(k)-v_i(k)\|^2  + \frac{1}{4c} \dist^2(v_i(k),\Xc).
\end{align}
By the definition of $x_i(k)$, we have $x_i(k) \in \mathcal{X}_i^{\Omega_i(k)}$, which implies
\[
\Es[\|v_i(k)-x_i(k)\|\mid \Fc_{k-1},I_k,J_k] \geq \Es[\text{dist}(v_i(k),\mathcal{X}_i^{\Omega_i(k)})\mid \Fc_{k-1},I_k,J_k].
\]
By Assumption \ref{assume:c} it follows
\[
\text{dist}^2(v_i(k),\mathcal{X}) \leq c\mathsf{E}\left[ \text{dist}^2(v_i(k),\mathcal{X}_i^{\Omega_i(k)})\mid\Fc_{k-1},I_k,J_k\right] \text{ for all } i.
\]
Therefore, the sum of the last two terms in Eq.~\eqref{eq:remove} is negative and by dropping that term,
we obtain the following relation w.p.1 for all $k\ge 1$ and $i \in \{I_k,J_k\}$,
\begin{align*}
\EXP{\|x_i(k)-x\|^2 \mid \Fc_{k-1},I_k,J_k}
\leq{}& (1-\rho_i) \|v_i(k)-x\|^2
- 2\alpha_i\langle \nabla f_i(x), z_i(k)-x\rangle
+ 8(1+c)\a_i^2G_f^2, 
\end{align*}
where $\rho_i=\sigma_i\alpha_i-4(2+c)\a_i^2L_i^2$.\quad$\blacksquare$

\subsection{Proof of Lemma \ref{lem:disagree2}\label{app:lemdisagree2}}
In the proof of Lemma~ \ref{lem:disagree2},
we use the following result (see Lemma 3.1 in \cite{Ram2010} for its proof).
\begin{lemma}\label{lem:scalar}
If $\lim_{k\rightarrow \infty} \g(k) = \g$ and $0<\beta<1$, then $\lim_{k \rightarrow \infty} \sum_{\ell=0}^k \beta^{k-\ell}\g(\ell) = \frac{\g}{1-\beta}$.
\end{lemma}
In addition, we also use an asymptotic upper bound
for the distance between the iterates
$x_i(k)$ and the set $\Xc$, which is given in the following lemma.
\begin{lemma}\label{lem:v}
Let Assumptions~\ref{assume:f}-\ref{assume:s} hold,
where the stepsizes $\a_i$ satisfy Assumption~\ref{assume:s}(a).
Then, for the iterates $x_i(k)$ of the method, we have
\[\limsup_{k \rightarrow \infty} \sum_{i=1}^m\Es[\dist^2(x_i(k),\Xc)]
\le\frac{8(1+c)m}{\min_i\{\g_i\rho_i\} }\bar{\g}\bar{\a}^2G_f^2,\]
where $\rho_i=\sigma_i\alpha_i-4(2+c)\a_i^2L_i^2.$
\end{lemma}
\begin{IEEEproof}
We use Lemma \ref{lem:base} with $x=\Ps_{\Xc}[v_i(k)]$, so that
w.p.1 for all $k\ge1$ and $i \in \{I_k,J_k\}$,
\begin{align*}
\EXP{\|x_i(k)-z_i(k)\|^2 \mid \Fc_{k-1},I_k,J_k}
\leq{}& (1-\rho_i) \dist^2(v_i(k),\Xc)
+ 8(1+c)\a_i^2G_f^2,
\end{align*}
with $\rho_i=\sigma_i\alpha_i-4(2+c)\a_i^2L_i^2$.
We note that
$\dist(x_i(k),\Xc)\le \|x_i(k)-z_i(k)\|$. Thus, we have w.p.1
for $i \in \{I_k,J_k\}$ and $k \ge 1$,
\begin{equation*}
\Es[\dist^2(x_i(k),\Xc)\mid \Fc_{k-1},I_k,J_k]
\leq \left(1-\rho_i\right)\dist^2(v_i(k),\Xc) + 8(1+c)\a_i^2G_f^2.
\end{equation*}
The preceding relation holds with probability $\g_i$ and, otherwise,
$x_i(k) = v_i(k)$ with probability $1-\g_i$. Thus, w.p.1 for all $k \ge 1$ and $i \in V,$
\begin{align*}
\Es[\dist^2(x_i(k),\Xc){}&\mid \Fc_{k-1}]
\leq (1-\g_i\rho_i)\,\EXP{\dist^2(v_i(k),\Xc) \mid\Fc_{k-1}}+ 8(1+c)\g_i\a_i^2G_f^2.
\end{align*}
By summing over $i$ and using relation (\ref{eqn:ds-dist}), we obtain
\begin{equation*}
\sum_{i=1}^m \Es[\dist^2(x_i(k),\Xc)\mid \Fc_{k-1}]
\leq \left(1-\min_i\{\g_i\rho_i\}\right)\sum_{j=1}^m \dist^2(x_j(k-1),\Xc)
+ 8(1+c)m \bar\g\bar\a^2G_f^2,
\end{equation*}
where $\bar{\g} = \max_i \g_i$ and $\bar{\a} = \max_i \a_i$.
Note that when $\rho_i\in (0,1)$ for all $i$, then we also have
$\g_i\rho_i\in (0,1)$ since $\g_i\in (0,1)$ for all $i$.
Taking the total expectation and, then, applying Lemma~\ref{lem:scalar} we obtain the desired
relation.
\end{IEEEproof}

\noindent
{\bf Proof of Lemma~ \ref{lem:disagree2}.}
We consider coordinate-wise relations similar to the proof of Lemma~\ref{lem:disagree}.
Since the matrices $W(k)$ are doubly stochastic for all $k\ge 1$, from relation (\ref{eqn:rec}) with $\|M\|=1$ and H\"older's inequality, we obtain
\begin{equation}\label{eqn:holder}
\sum_{\ell=1}^d \Es[\|y_{\ell}(k)-[\bar{x}(k)]_{\ell}\1b\|^2]\le
\left(\sqrt{\sum_{\ell=1}^d \Es[\|D_k(y_{\ell}(k-1)-[\bar{x}(k-1)]_{\ell}\1b)\|^2]} + \sqrt{\sum_{\ell=1}^d \Es[\|\delta_{\ell}(k)\|^2]}\right)^2,
\end{equation}
where $[\bar{x}(k)]_{\ell} = \frac{1}{m}\1b^Ty_{\ell}(k)$, $D_k = W(k) - \frac{1}{m}\1b\1b^T$,
and $\lambda < 1$. From relation~\eqref{eq:ly}, we know that
\begin{equation}\label{eqn:eqrec3}
\sum_{\ell=1}^d \Es\left[\|D_k(y_{\ell}(k-1)-[\bar{x}(k-1)]_{\ell}\1b)\|^2\right] \le \lambda \sum_{\ell=1}^d \Es\left[\|y_{\ell}(k-1)-[\bar{x}(k-1)]_{\ell}\1b\|^2\right].
\end{equation}
The second term in (\ref{eqn:holder}) is evaluated similar to that in Lemma \ref{lem:disagree}.
Hence, for all $k \ge 1$ w.p.1,
\begin{align}\label{eqn:eqrec4}
\sqrt{\sum_{\ell=1}^d\Es\left[\|\delta_{\ell}(k)\|^2\right]}  {}\le \beta_k,
\end{align}
where
$\beta_k=\sqrt{(4+ 4\bar{\a}^2\bar L^2)\sum_{i=1}^m \Es\left[\dist^2(v_i(k),\Xc)\right]+ 4m\bar{\a}^2G_f^2},$
$\bar{\a} = \max_i \a_i$ and $\bar L=\max_i L_i.$
Letting $u_k = \sqrt{\sum_{\ell=1}^d \Es[\|y_{\ell}(k)-[\bar{x}(k)]_{\ell}\1b\|^2]}$ in~(\ref{eqn:holder}) and
using relations (\ref{eqn:eqrec3}) and (\ref{eqn:eqrec4}), we have
\[
u_k \le \sqrt{\lambda}
u_{k-1} + \beta_k\qquad\hbox{for all $k \ge 1$}.
\]
Since $\lambda <1$, by Lemma \ref{lem:scalar} we have
$\limsup_{k \rightarrow \infty} u_k \le \limsup_{k \to\infty}\beta_k/(1-\sqrt{\lambda}),$
implying that
\begin{equation}\label{eq:l1}
\limsup_{k \rightarrow \infty} u_k^2 \le \frac{1}{(1-\sqrt{\lambda})^2} \limsup_{k \to\infty}\beta_k^2.
\end{equation}
By relation~\eqref{eqn:ds-dist} it follows that
\begin{equation}\label{eq:l2}
\limsup_{k\to\infty}\b^2_k\le
(4+ 4\bar{\a}^2\bar L^2) \limsup_{k\to\infty} \sum_{j=1}^m \Es[\dist^2(x_j(k-1),\Xc)]
+4m\bar{\a}^2G_f^2.\end{equation}
Finally, by the definition of $y_{\ell}(k)$, we have $u_k^2 = \sum_{i=1}^m \Es[\|x_i(k)-\bar{x}(k)\|^2]$.
The desired relation follows from Eqs.~\eqref{eq:l1} and~\eqref{eq:l2} and Lemma \ref{lem:v}.
\quad $\blacksquare$
\end{document}